%% file: paperV2.tex
\documentclass[a4wide,twoside,12pt]{article}

\usepackage{a4wide}
\usepackage{graphicx}
\usepackage{amstext, amsmath, amssymb, amsfonts, amsbsy}

\usepackage{latexsym}
\usepackage{booktabs}
\usepackage{url}
\usepackage{import}
\usepackage{fancyvrb}
\usepackage{stmaryrd}
\usepackage{pdfsync}
\usepackage{algorithm}
\usepackage{tikz}
\usepackage{pgfplots}
\usepackage{pgfplotstable}
\usepackage{subfig}
\usepackage{tabularx}
\usepackage{snapshot}

\usepackage{array}
\newcolumntype{L}{@{}>{\kern\tabcolsep}l<{\kern\tabcolsep}}
\usepackage{colortbl}
\usepackage{xcolor}

\usepackage{algpseudocode}

\usepackage[shadow]{todonotes}
\usepackage[numbers]{natbib}

\numberwithin{equation}{section}
\numberwithin{figure}{section}
\numberwithin{table}{section}

\usepackage{amsthm}

\input{notation.tex}

\title{\bf A CutFEM method for two-phase flow problems}

\author{Susanne Claus, Pierre Kerfriden
\\ \\
Cardiff University, School of Engineering, \\ The Parade, CF243AA Cardiff, United Kingdom \\
}

\begin{document}
\maketitle
\renewcommand{\thefootnote}{\fnsymbol{footnote}}


\renewcommand{\thefootnote}{\arabic{footnote}}

\begin{abstract}
  \noindent 
In this article, we present a cut finite element method for two-phase Navier-Stokes flows. The main feature of the method is the formulation of a unified continuous interior penalty stabilisation approach for, on the one hand, stabilising advection and the pressure-velocity coupling and, on the other hand, stabilising the cut region. The accuracy of the algorithm is enhanced by the development of extended fictitious domains to guarantee a well defined velocity from previous time steps in the current geometry. Finally, the robustness of the moving-interface algorithm is further improved by the introduction of a curvature smoothing technique that reduces spurious velocities. The algorithm is shown to perform remarkably well for low capillary number flows, and is a first step towards flexible and robust CutFEM algorithms for the simulation of microfluidic devices.
\end{abstract}


\section{Introduction}
The development of general-purpose, robust and accurate finite element schemes for the simulation of two-phase immiscible fluids is hindered by several key scientific challenges. Firstly, the interface between the two fluids needs to be tracked whilst allowing topological changes such as droplet merging or breakup. Secondly, a contrast in fluid densities and viscosities causes a velocity kink at the interface. Similarly, surface tension forces cause the pressure field to jump at the interface. Discontinuities of these types are not easily represented in a finite element context, where fields are assumed to be piecewise smooth by construction. Finally, these difficulties, which are specific to multi-phase flows, have to be solved in the context of general fluid dynamics, where advection and pressure-velocity coupling are well-known sources of numerical instabilities. 

Two main types of numerical approaches can be distinguished to represent the interface between two fluids: sharp interface approaches and diffusive interface approaches. In diffusive interface approaches (e.g. phase field methods~\cite{Anderson1998,Jacqmin1999}, conservative level set methods~\cite{Olsson2005}, level set methods~\cite{Sussman1994}), a smoothing region around the interface is introduced, which blends viscosities and densities smoothly from one phase to the other. In this context, surface tension forces are usually applied over the entire smoothing region (e.g. continuum surface force approach \cite{Brackbill1992}). Smoothing the interface makes the solver more robust, in particular when surface tension is strong. Moreover, topological changes are handled in a straightforward manner. However, the most stringent drawback of this type of approach is the need for a high number of elements to resolve the smoothing region. \\
In contrast, sharp interface methods, which are of interest in this article, do not introduce such a smoothing region and hence less elements are required to resolve the interface region. 
The most well-known sharp interface method is the ALE method \cite{DoneaALE}, in which the elements of the finite element mesh are deformed to match the fluid interface. This mesh moving approach works remarkably well for small deformations but requires complex re-meshing procedures for large deformations and topological changes. Alternatively, an emerging family of methods allow the sharp interface
to cut through the elements, and to move through a fixed background grid without re-meshing. In this case, the immersed interface may be tracked using either Lagrangian markers (e.g. immersed boundary method~\cite{Mittal2005}, front-tracking~\cite{Unverdi1992}) or functional representations (e.g. level-set-method \cite{Sussman1994}). Importantly, elements which are intersected by the interface need to be enriched in order to capture jumps and kinks in the solution. Failure to do so significantly deteriorates convergence with mesh refinement and introduces significant spurious velocity modes for flows involving surface tension. Finally, the most complex aspect of this family of methods is their intrinsic need for tailored stabilisation. XFEM and CutFEM are two popular frameworks proposed to systematise the process of embedding discontinuities in finite element solvers, which may be used to build two-phase flow solvers. XFEM enriches the finite element shape functions by the Partition-of-Unity method \cite{Melenk1996} and was originally introduced by \cite{Moees1999} for crack propagation problems. XFEM was first developed for solid mechanics and was later extended to two-phase flow problems \cite{Chessa2003,Gross2006, Fries2009,Sauerland2013}. Stabilisation in XFEM is usually achieved through a combination of pre-conditioning \cite{Kirchhart2016, Gros2016, Lehrenfeld2017}, deactivation of degrees of freedom \cite{Fries2006} or more recently the use of ghost penalty stabilisation \cite{Schott2015, Gros2016}. The CutFEM approach \cite{BurmanClausHansboEtAl2014} uses overlapping fictitious domains in combination with ghost penalty stabilisation \cite{Burman2010} to enrich and stabilise the solution. CutFEM has been developed for, e.g. Stokes interface problems \cite{Hansbo2014, MassingLarsonLoggEtAl2013}, for the Oseen problem \cite{Winter2018} and fluid structure interaction problems \cite{Massing2015,Burman2014}. An alternative to the two families of methods described previously was proposed in \cite{Heimann2013}, where an unfitted discontinuous Galerkin approach was successfully developed to solve two-phase flow problems.

In this article, we develop a new CutFEM solver for multiphase flows. The main novel and appealing feature of the proposed framework is the unification of the stabilisation processes. More precisely, we will use continuous interior penalty (CIP) to stabilise the unfitted interface discretisation, which is known as ghost-penalty in the context of CutFEM, and we will also use CIP to stabilise advection and the pressure-velocity coupling. In addition to the development, implementation and validation of this unified framework, we  propose a series of novel tools to improve the capabilities and performances  of the CutFEM framework when addressing two-phase flows problems. These are:
(i) the definition of the interface tracking problem using extended fictitious subdomains. This is to guarantee that velocities from the past time step are well-defined in the current geometry; and
(ii) a curvature smoothing for the surface tension force which reduces the amplitude of spurious interface velocity modes.
Consistently with our previous CutFEM developments \cite{BurmanClausHansboEtAl2014,Claus2018,Claus2018a}, the different phases are described implicitly using the zero isoline of a level-set function. Finite element enrichment is obtained
through the application of an overlapping domain decomposition technique that allows the pressure to jump  and the velocity to kink across the fluid interface. Time integration is performed using an implicit Euler scheme, whilst the velocity and pressure fields are spatially discretised using linear continuous finite element spaces. We validate the developments in several ways. Firstly, we investigate the performance of the two-phase flow solver on the rising bubble benchmark. We show that our results agree with that of existing schemes \cite{Hysing2009}. Then, we study and discuss the effect of the parameters of the proposed unified stabilisation scheme on the shape of the rising bubble. Finally, we solve the problem of a droplet in a microfluidic 5:1:5 contraction and expansion flow, which allows us to prove the robustness of the CutFEM scheme for a surface tension dominated flow.

The paper is organised as follows: In Section 2, we introduce the strong formulation of the two-phase Navier-Stokes problem together with the level set advection problem followed by their weak formulation. In Section 3, we introduce the extended fictitious domains, describe the discretisation in time  and introduce the stabilised cut finite element scheme for two-phase Navier-Stokes flows. We then discuss the stabilised curvature projection before summarising the scheme in Section 4. In Section 5, we present numerical results for the rising bubble benchmark and a droplet in a 5:1:5 contraction and expansion problem.

\clearpage
\section{The Two-Phase Navier-Stokes Problem}
Let $\Omega$ be a domain in $\mathbb{R}^d$ ($d=2$ or $3$) with a
polyhedral boundary $\partial \Omega$. We assume that two immiscible
incompressible Newtonian fluids occupy the time-dependent subdomains $\Omega_i = \Omega_i(t)
\subset \Omega$, $i=1,2$, with $\Omega = \Omega_1(t) \cup \Omega_2(t)$ and
that $\Gamma = \Gamma(t)$ denotes the smooth interface between them.  
\paragraph{Two-phase Navier-Stokes Problem}
The two-phase flow system is described by the following Navier-Stokes problem in the time interval $t \in (0,T]$: Find the velocity
$\vel: \Omega \times (0,T] \rightarrow \mathbb{R}^d$ and the pressure $p: \Omega \times (0,T]
\rightarrow \mathbb{R}$, such that
\begin{align}
 \rho \frac{\partial \vel}{\partial t} + \rho(\vel \cdot \nabla) \vel  - \nabla \cdot \stress &= \rho \mathbf{\bfg}  &\mbox{ in } \Omega, \nonumber \\ 
\stress  &= 2 \eta \epsilon(\vel) - p \mathbf{I} &\mbox{ in } \Omega, \nonumber \\
\nabla \cdot \vel &= 0 &\mbox{ in } \Omega .
\end{align}
Here, $\stress$ is the Cauchy stress tensor, $\epsilon(\vel) = \frac{1}{2} \left(\nabla \vel + \nabla \vel^T\right)$ is
the rate of deformation tensor, $\eta|_{\Omega_i} = \eta_i$, $i=1,2$,
is the piecewise constant fluid viscosity on $\Omega_i$, $\rho|_{\Omega_i} = \rho_i$, $i=1,2$, is the piecewise constant density on $\Omega_i$ and $\bfg \in
[L^2(\Omega)]^d$  is a known force (e.g. gravity).
Additionally, we have given initial conditions for the velocity $\vel(\bfx,0) = \vel_0$ and boundary conditions on $\partial \Omega$ which we here assume to be of Dirichlet type for simplicity, i.e. 
\begin{equation}
\vel = \vel_D \mbox{ on } \partial \Omega
\end{equation}
with given function $\vel_D \in  [H^{1/2}(\partial \Omega)]^d $. At the interface between the two fluids, we enforce a zero jump
condition on the velocity, $\jump{\vel} = (\vel_1 - \vel_2)|_{\Gamma}
= 0$, $\vel_i=\vel|_{\Omega_i}$, $i=1,2$, and the jump of the traction
vector over the interface is given by the surface tension, i.e.
\begin{align}
\jump{\vel} &= 0 &\mbox{ on } \Gamma(t), \nonumber \\
\jump{\stress \cdot \bfn_{\Gamma}} &= -\gamma \kappa \bfn_{\Gamma} &\mbox{ on } \Gamma(t).
\label{equ: surface tension}
\end{align}
Here,
$\gamma$ is the surface tension coefficient, $\kappa$ is the curvature
of the interface and $\bfn_{\Gamma}$ is the outward unit normal on the
interface pointing from $\Omega_1(t)$ to $\Omega_2(t)$. \paragraph{Interface Evolution}
To capture the evolution of the interface $\Gamma(t)$, we will use a level set function. The level set function is a scalar function on $\Omega$, such that $\phi(\bfx,t)<0$ for  $\bfx \in \Omega_1(t)$, $\phi(\bfx,t)>0$ for  $\bfx \in \Omega_2(t)$ and $\phi(\bfx,t)=0$ for  $\bfx \in \Gamma(t)$. The level set function will be moved with the fluid velocity $\vel$ by solving the following advection problem: Find the level set function $\phi: \Omega \times (0,T] \rightarrow \mathbb{R}$, such that 
\begin{equation}
\frac{\partial{\phi}}{\partial t} + \vel \cdot \nabla \phi = 0
\end{equation}
with given initial condition $\phi(\bfx,0) = \phi_0$ and inflow boundary conditions $\phi = \phi_{in}$ at the inflow boundary $\partial \Omega_{in}(t):= \{\bfx \in \partial \Omega(t): \vel(\bfx,t) \cdot \bfn <0\}$. Here, $\bfn$ is the outward normal to $\Omega$.

\subsection{Weak Formulation}
Let us denote  
\begin{equation}
\begin{aligned}
\left( u, v \right)_{\Omega_i} = \int_{\Omega_i} u_i v_i \, dx \\
\left( \jump{u}, \jump{v} \right)_{\Gamma} = \int_{\Gamma} \jump{u} \jump{v} \, ds
\end{aligned}
\end{equation}
and
\begin{equation}
\left( u, v \right)_{\Omega} = \sum_{i} \left( u_i, v_i \right)_{\Omega_i}, \quad i=1,2.
\end{equation}
And let us introduce the following Sobolev space
\begin{equation}
\begin{aligned}
L^2_0(\Omega) &= \left\{ p \in L^2(\Omega): \int_{\Omega}  p = 0\right\}. \\
\end{aligned}
\end{equation}
\paragraph{Weak Formulation of the Two-Phase Navier-Stokes Problem}
We seek the solution to the two-phase Navier-Stokes problem for $t \in (0,T]$, for $U(t):=(\vel(t),p(t)) \in [H^1(\Omega)]^d \times L^2_0(\Omega)$, such that 
\begin{equation}
A(U(t),V) = L(V)
\end{equation}
for all $V=(\velt,q) \in [H^1(\Omega)]^d \times L^2_0(\Omega)$. \\
Here, 
\begin{equation}
\begin{aligned}
A(U(t),V) &= \left(\rho \frac{\partial \vel(t)}{\partial t}, \velt\right)_{\Omega} + \mathcal{N}(\vel(t),\vel(t),\velt)  \\
&\quad+ a(\vel(t),\velt) + b(p(t),\velt) -  b(q,\vel(t)) + a_n(\vel(t),\velt),  \\
L(V) &= (\rho \bfg,  \velt)_{\Omega} - (\gamma \kappa \bfn_{\Gamma}, \mean{\velt})_{\Gamma} + l_n(\velt),
\end{aligned}
\end{equation}
consist of a non-linear convection term 
\begin{equation}
\begin{aligned}
\mathcal{N}(\vel(t),\vel(t),\velt) &= \left( \rho (\vel(t) \cdot \nabla) \vel(t), \velt\right)_{\Omega},
\end{aligned}
\end{equation}
terms that arise from the integration by parts of $(\nabla \cdot \stress,\velt)_{\Omega}$
\begin{equation}
\begin{aligned}
a(\vel(t),\velt) &= \left( 2 \eta \epsilon(\vel(t)),\epsilon(\velt) \right)_{\Omega}  -  \left( 2 \eta \epsilon(\vel(t)) \cdot \bfn, \velt \right)_{\partial \Omega}  -  \left( \avg{2 \eta \epsilon(\vel(t)) \cdot \bfn_{\Gamma}},\jump{\velt} \right)_{\Gamma} , \\ 
b(p(t),\velt) &= - (p(t), \nabla \cdot \velt)_{\Omega} + ( p(t) \cdot \bfn, \velt)_{\partial \Omega} + (\avg{p(t) \cdot \bfn_{\Gamma}}, \jump{\velt})_{\Gamma}
\end{aligned}
\end{equation}
and terms that enforce the interface condition $\jump{\vel} = \mathbf{0}$ and the boundary condition $\vel = \vel_D$ weakly using Nitsche's method
 \begin{equation}
\begin{aligned}
a_n(\vel,\velt) &=   -  \left( 2 \eta \epsilon(\velt) \cdot \bfn ,\vel \right)_{\partial \Omega}  + \lambda_{\partial \Omega}(h) \, \eta \left( \vel, \velt \right)_{\partial \Omega}   \\
&\phantom{=}  -  \left( \avg{2 \eta \epsilon(\velt) \cdot \bfn_{\Gamma}},\jump{\vel} \right)_{\Gamma} + \lambda_{\Gamma}(h) \, \avg{\eta} \left( \jump{\vel},\jump{\velt} \right)_{\Gamma}, \\
l_n(\velt) & =  -  \left( 2 \eta \epsilon(\velt) \cdot \bfn, \vel_D \right)_{\partial \Omega}  -  ( q \cdot \bfn, \vel_D)_{\partial \Omega} +  \lambda_{\partial \Omega}(h) \, \eta \left( \vel_D, \velt \right)_{\partial \Omega}. 
\end{aligned}
\end{equation}
Here, $\lambda_{\Gamma}(h) \in \mathbb{R}$ and  $\lambda_{\partial \Omega}(h) \in \mathbb{R}$  are positive penalty
parameters that have to be chosen large enough to guarantee the
coercivity of the weak formulation and $\bfn$ is the outward normal on $\partial \Omega$. 

Note that the interface terms in the equations above have been obtained through integration by parts of $(\nabla \cdot \stress,\velt)_{\Omega}$, summing over $\Omega_i$ and then using the relationship 
\begin{equation}
\jump{ab} = \jump{b}\avg{a}+\mean{b}\jump{a}
\end{equation}
to reformulate them in terms of jumps and averages. The averages are weighted and defined as 
\begin{align}
\avg{\bfv} := w_1 \bfv_1 + w_2 \bfv_2,  \\
\mean{\bfv} := w_2 \bfv_1 + w_1 \bfv_2,
\end{align}
where $w_1$ and $w_2$ are weights to be chosen such that  $w_1+w_2 =1, 0 \leq w_i \leq 1, i =1,2$. 
These weights need to be chosen carefully for problems involving a high contrast in material parameters (see e.g. \cite{Burman2011} for advection diffusion problems). To ensure robustness with respect to a high contrast in the fluid viscosities as well as fluid densities, we choose the weights as the harmonic average over the kinematic viscosity $\nu_i:=\frac{\eta_i}{\rho_i}$, $i=1,2$, i.e.
\begin{equation}
\begin{aligned}
w_1 := \frac{\nu_2}{\nu_1+ \nu_2}, \\ 
w_2 := \frac{\nu_1}{\nu_1+ \nu_2}.
\end{aligned}
\end{equation}
\paragraph{Weak Formulation of the Level Set Advection}
The weak formulation for the level set advection problem reads: For $t \in (0,T]$, find $\phi(t) \in H^1(\Omega)$, such that 
\begin{equation}
\left( \frac{\partial \phi(t)}{\partial t}, \psi \right)_{\Omega} + \left( \vel(t) \cdot \nabla \phi(t), \psi \right)_{\Omega} = 0
\end{equation}
for all $\psi \in H^1(\Omega)$ with suitable initial and boundary conditions.

\section{Stabilised Cut Finite Element Formulation}
\subsection{Mesh Quantities} 
\begin{figure}[htbp]
\centering
\includegraphics[width=.9\textwidth]{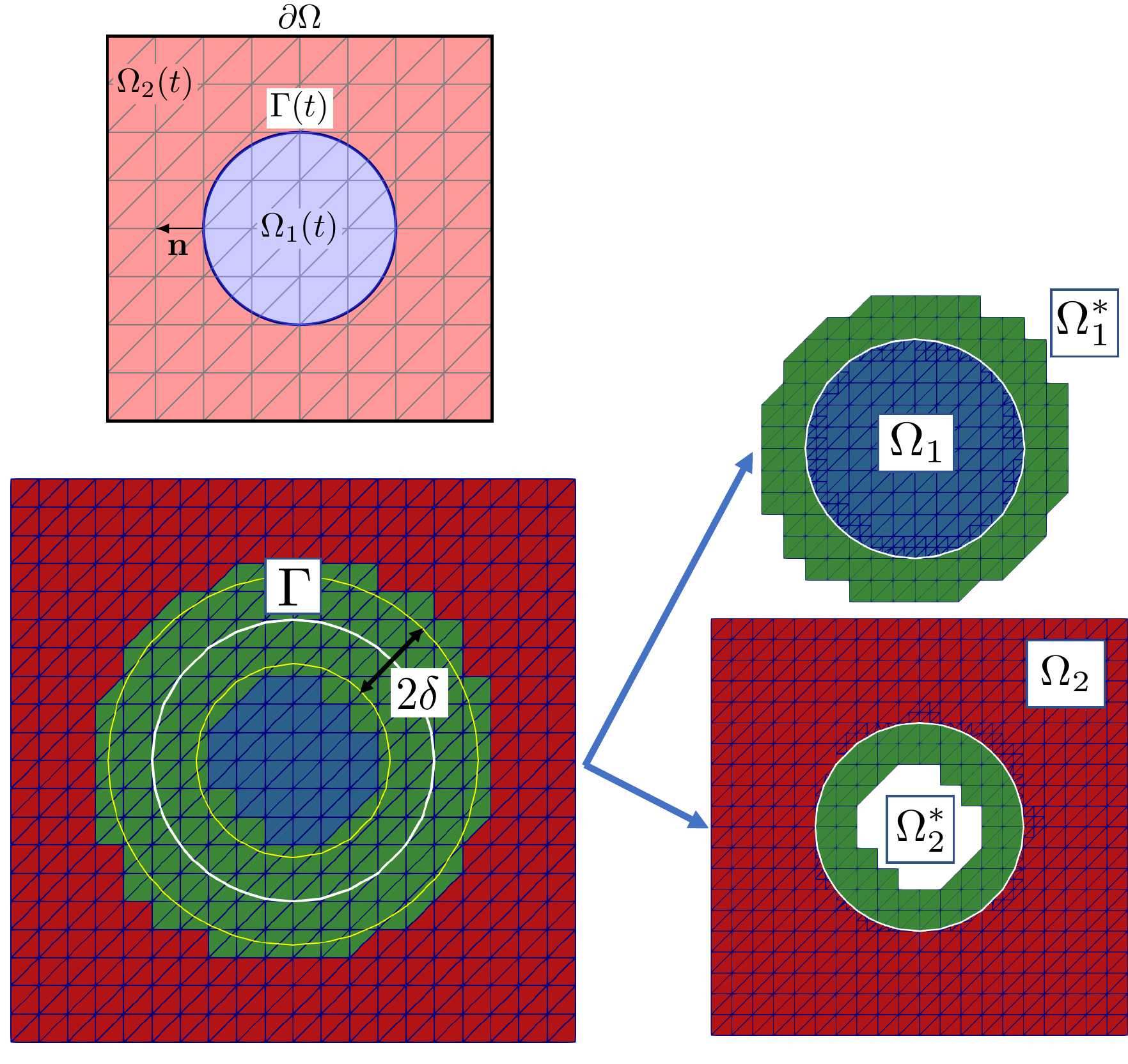}
 \caption{Navier-Stokes multi-domain discretisation. Cells with a distance smaller than $\delta$ to the interface are doubled (elements shaded in green). The interior domain $\Omega_1$ gets extended to the fictitious domain $\Omega_1^*$ and the exterior domain, $\Omega_2$,  gets extended to the fictitious domain $\Omega_2^*$  resulting in overlapping fictitious domains (green shaded elements).}
 \label{fig: domains}
\end{figure}

\begin{figure}[htbp]
\centering
\subfloat[]{\includegraphics[width=.9\textwidth]{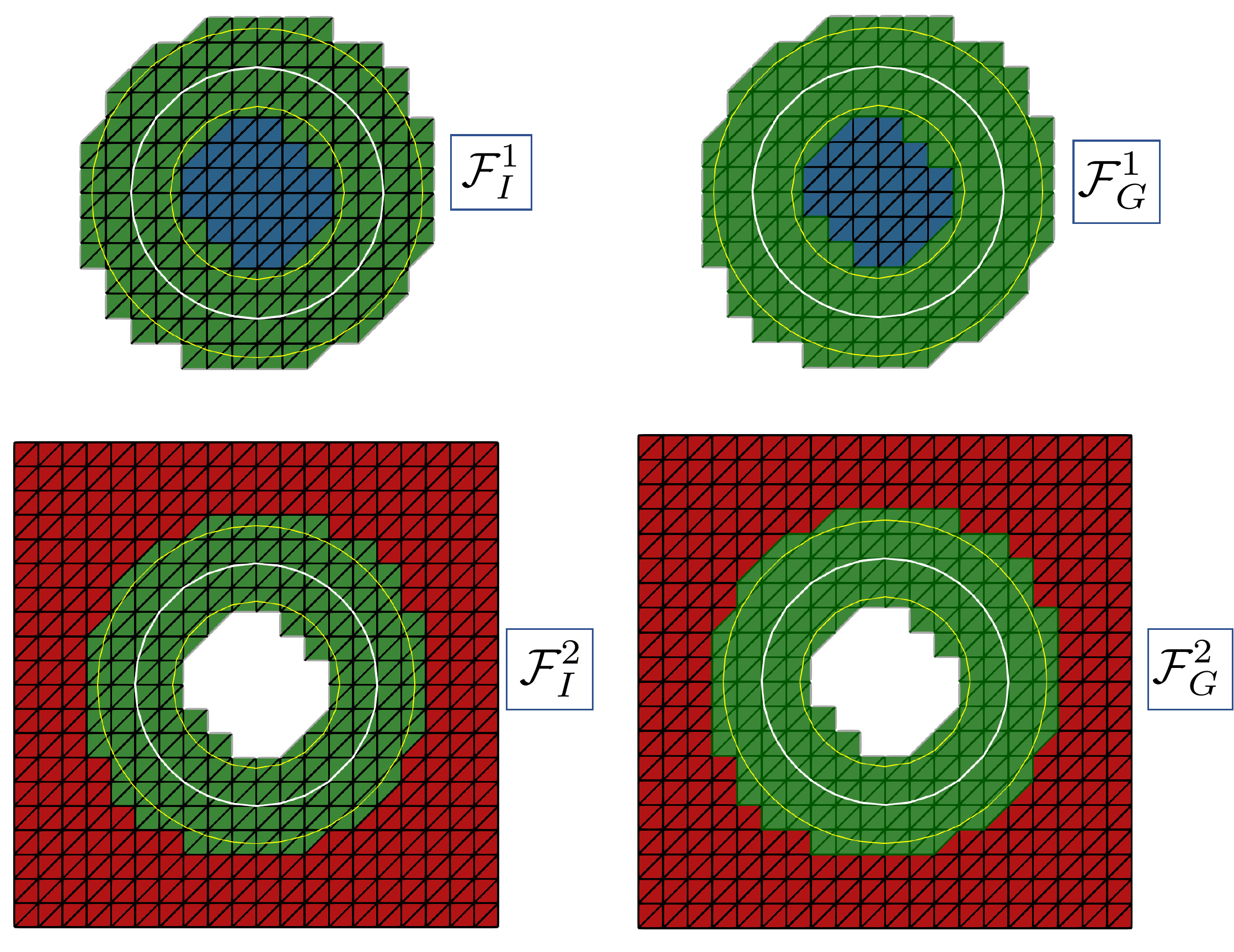}}\\
\subfloat[]{\includegraphics[width=.4\textwidth]{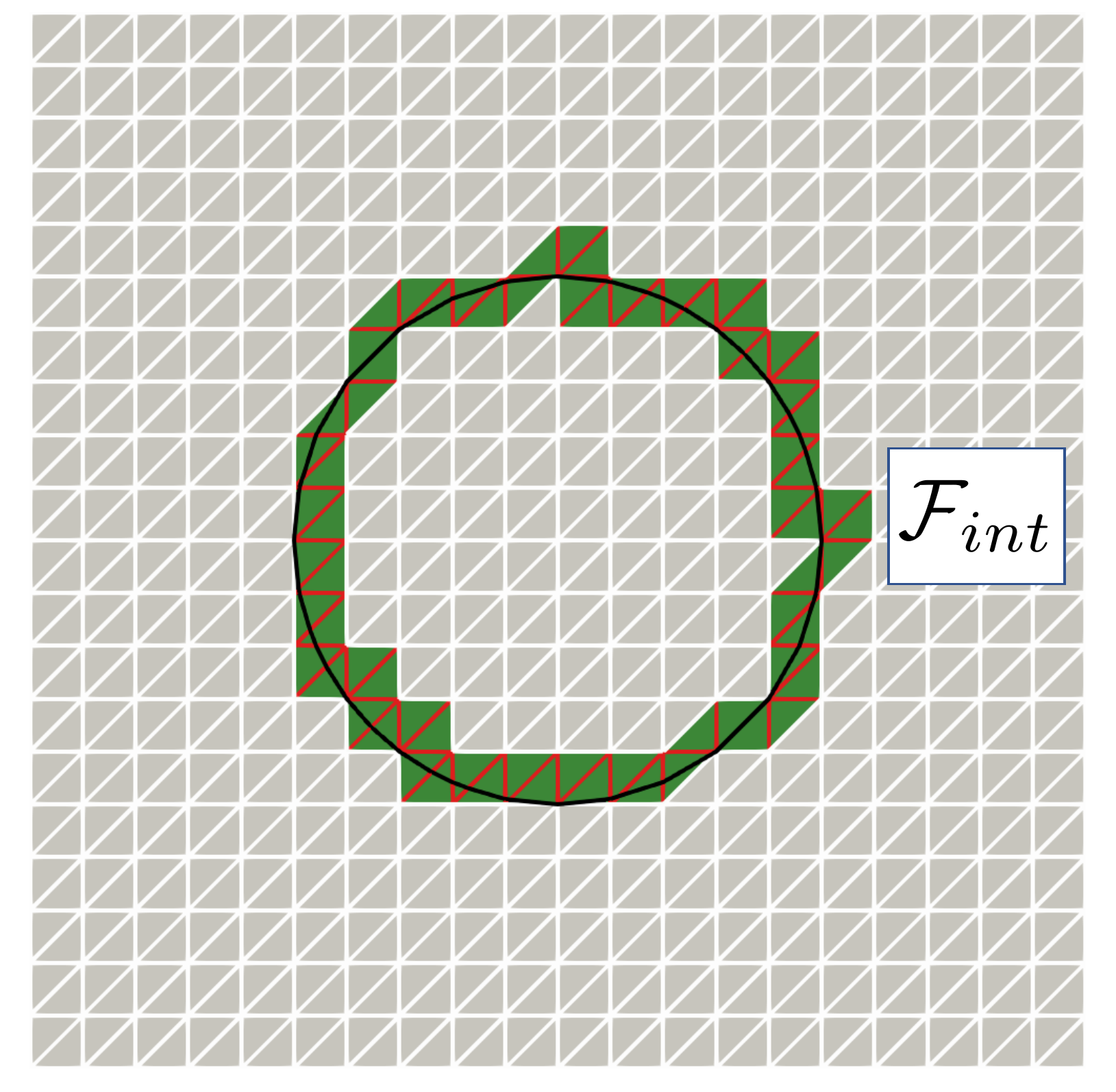}}
 \caption{Illustration of different facet sets in two-phase Navier-Stokes problem. The top figure
 illustrates the definition of interior $\mathcal{F}_{I}^i$ (black edges) and ghost
 penalty facets  $\mathcal{F}_{G}^i$ (dark green edges) and the bottom figure displays the facets which are intersected by the interface, $\mathcal{F}_{int}$ (red edges).}
 \label{fig: facets}
\end{figure}

In this section, we define mesh quantities that will be used in our finite element discretisation.  We discretise the domain $\Omega$, which will stay fixed in time, using a quasi-uniform triangulation, $\left\{ \mathcal{T}_h \right \}_{0 < h \leq 1}$, independent of the location of the interface $\Gamma(t)$ but coinciding with the outer boundary $\partial \Omega$. We consider the evolution of $\Gamma(t)$ through $\Omega$ in the time interval $t \in [0,T]$, which we decompose into  $n_t$ equal steps, i.e. $t_n= \Delta t \, n$ with $\Delta t = T/n_t$.  For each $t_n$, we define a range of sets of active elements and facets (i.e. faces in 3D and edges in 2D) in $\mathcal{T}_h $ associated with $\Gamma^n = \Gamma(t_n)$ and $\Omega^n_i = \Omega_i(t_n)$. \\
Firstly, let us consider the $\delta$ neighbourhood of $\Gamma^n$ defined as 
\begin{equation}
\mathcal{O}_{\delta}^n = \{ \bfx \in \mathbb{R}^d : dist(\bfx,\Gamma^n) \leq \delta\} 
\end{equation}
and the extended domains 
\begin{equation}
\Omega_{i, \delta}^n = \Omega_i^n \cup \mathcal{O}_{\delta}^n, i =1,2.
\end{equation}
Secondly, let us define submeshes of $\mathcal{T}_h$, associated with these extended domains, consisting of all elements in $\mathcal{T}_h$ that have at least a small part in domain $\Omega_{i,\delta}^n$ as
\begin{equation}
\mathcal{T}_{h,i}(t_n) = \{ K \in \mesh : K \cap \Omega_{i,\delta}^n \neq \emptyset\}, i=1,2.
\end{equation}
We call $\mathcal{T}_{h,i}(t)$ the \emph{active mesh} of $\Omega_i(t)$. Note that these active meshes do not conform to $\Gamma(t)$ as shown in Figure~\ref{fig: domains}. We denote the union of all elements in $\mathcal{T}_{h,i}(t_n)$ as 
\begin{equation}
\Omega^*_i(t_n) = \bigcup_{K \in \mathcal{T}_{h,i}(t_n)} K, i=1,2.
\end{equation}
$\Omega^*_i(t)$ is called the \emph{fictitious domain} of $\Omega_i(t)$ (see Figure~\ref{fig: domains}).
Furthermore, consider the set of all elements with part in the $\delta$ neighbourhood, denoted by
\begin{equation}
G_h^n = \{ K \in \mesh : K \cap \mathcal{O}_{\delta}^n \neq \emptyset\}.
\end{equation}
We will apply stabilisation terms to a subset of facets (i.e. faces in 3D and edges in 2D) in these domains. In particular, for each active mesh at time $t_n$, we will apply stabilisation terms to so-called \emph{interior facets} and \emph{ghost penalty facets} of $\mathcal{T}_{h,i}(t_n)$. The set of interior facets, $\mathcal{F}_I^{i}(t_n)$, contains all facets that are shared by two elements in $\mathcal{T}_{h,i}(t_n), i =1,2$ (see Figure~\ref{fig: facets}). The set of intersected facets, $\mathcal{F}_{int}(t_n)$, contains all facets intersected by the interface. And ghost penalty facets are all facets in the set 
\begin{equation}
\mathcal{F}^i_{G}(t_n) = \{ F \in  \mathcal{F}_I^i(t_n): K^+ \in G_h^n \mbox{ or } K^- \in G_h^n \}. 
\end{equation}
Here, $K^+$ and $K^-$ are the two elements sharing the interior
face $F \in \mathcal{F}_I^i(t_n)$ (see Figure~\ref{fig: facets}).  \\
We place the following
assumptions on the background mesh $ \mathcal{T}_h$, to ensure that
the interface $\Gamma(t)$ is sufficiently resolved: For $t \in [0,T]$,
\begin{itemize}
\item G1:  $\Gamma(t)$ crosses element faces at most once.
\item G2: For each element $K$ intersected by $\Gamma(t)$, there exists a
  plane $S_K$ and a piecewise smooth parametrization $\Phi: S_K \cap K
  \rightarrow \Gamma(t) \cap K$.
\item G3: We assume that for each intersected element $K $,
  there are elements $K_i \in \Omega_i(t), i=1,2$, such that
  $\overline{K} \cap \overline{K_i} \neq \emptyset$.  This means each
  intersected element shares at least one vertex or facet with an
  element $K_1$ in $\Omega_1(t)$ and an element $K_2$ in $\Omega_2(t)$.
 \end{itemize}

\subsection{Discretisation in Time} 
As mentioned in the previous section, we decompose the time interval $[0,T]$ into $n_t$ equally spaced time steps. For each $t_n$, we decouple the two-phase Navier-Stokes problem from the level set advection by first solving the two-phase Navier-Stokes problem in $\Omega = \Omega^n_1 \cup \Omega^n_2$ and then use the obtained velocity to advect the level set function to obtain the updated domains $\Omega_i^{n+1}$. \\
We employ an implicit Euler scheme to discretise the two-phase Navier-Stokes problem in time which yields 
\begin{equation}
\begin{aligned}
 \rho \left( \frac{\vel^{n}-\vel^{n-1}}{\Delta t}, \velt\right)_{\Omega} + \mathcal{N}(\vel^{n-1},\vel^{n},\velt) 
+ a(\vel^{n},\velt) + a_n(\vel^{n},\velt) + b(p^{n},\velt) -  b(q,\vel^{n}) = L(V) 
\end{aligned}
\end{equation}
together with suitable initial conditions. 
Here, we have linearised the non-linear term using the approximation 
\begin{equation}
\mathcal{N}(\vel^{n},\vel^{n},\velt) \approx \mathcal{N}(\vel^{n-1},\vel^{n},\velt)  =  \rho \left( (\vel^{n-1} \cdot \nabla) \vel^{n}, \velt\right)_{\Omega}. 
\end{equation}
Note that all integrals above are evaluated in the domains $\Omega_i(t_n)$, i.e. 
\begin{equation}
\left ( \vel, \velt \right)_{\Omega} = \sum_{i=1}^2 \left ( \vel_i , \velt_i \right)_{\Omega_i^n}.  
\end{equation}
This poses the challenge to evaluate the velocity from the previous time step, $\vel^{n-1}$, which was obtained in $\Omega_i(t_{n-1})$ in the updated domains $\Omega_i(t_n)$. To achieve this, we define an extension of $\vel^{n-1}$ from $\Omega_i(t_{n-1})$ to the extended domains $\Omega_{i,\delta}^{n-1}$ using an extension operator
\begin{equation}
\mathcal{E}^{n-1}_i: H^1(\Omega_i^{n-1}) \rightarrow H^1(\Omega_{i,\delta}^{n-1}).  
\end{equation}
And we require the $\delta$-neighbourhood to be large enough such that 
\begin{equation}
\Omega_{i}^{n} \subset  \Omega_{i,\delta}^{n-1}. 
\end{equation}
This guarantees that the extensions of $\vel_i^{n-1}$, denoted by $\mathcal{E}^{n-1}_i \vel^{n-1}_i$, $i=1,2$, are well defined in $\Omega_i(t_n)$. We introduce the velocity of the previous time-step evaluated as its extension in $\Omega = \Omega_1^n \cup \Omega_2^n$ as 
\begin{equation}
\bfbeta := 
\begin{cases}
\mathcal{E}^{n-1}_1 \vel^{n-1}_1 & \mbox{for } \bfx \in \Omega_1^n, \\
\mathcal{E}^{n-1}_2 \vel^{n-1}_2 & \mbox{for } \bfx \in \Omega_2^n \backslash \Gamma.  
\end{cases}
\end{equation}
Using this extension definition, we reformulate the two-phase Navier-Stokes problem as 
\begin{equation}
\begin{aligned}
 \rho \left( \frac{\vel^{n}-\bfbeta}{\Delta t}, \velt\right)_{\Omega} + \mathcal{N}(\bfbeta,\vel^{n},\velt) 
+ a(\vel^{n},\velt) + a_n(\vel^{n},\velt) + b(p^{n},\velt) -  b(q,\vel^{n}) = L(V). 
\end{aligned}
\end{equation}
The level set advection problem is discretised in time using a $\theta$-scheme \cite{Gross2011} which yields
\begin{equation}
\left( \frac{\phi^{n}-\phi^{n-1}}{\Delta t}, \psi \right)_{\Omega} + \theta \left( \vel^n \cdot \nabla \phi^n, \psi \right)_{\Omega} + (1-\theta) \left( \bfbeta \cdot \nabla \phi^{n-1}, \psi \right)_{\Omega} = 0
\label{equ: level set time discretisation}
\end{equation}
with suitable initial and boundary conditions. In this contribution, we choose $\theta =0.5$ which yields a Crank-Nicolson scheme.

\subsection{Discretisation in Space}
\subsubsection{Domain and interface approximation}

\begin{figure}
\centering
\includegraphics[width=.5\textwidth]{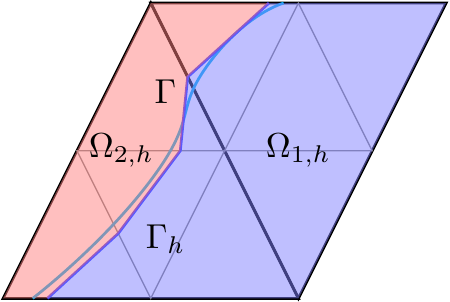}
\caption{Illustration of linear approximation of interface $\Gamma$ and domains $\Omega_i, i=1,2$, on the refined mesh $\mathcal{T}_{h/2}$ with respect to the mesh $\mesh$.}
\label{fig: refined mesh approximation}
\end{figure}

\paragraph{Domain approximation} 
As previously mentioned, we describe our domains $\Omega_1(t)$ and $\Omega_2(t)$ using a level set function. This level set function is discretised using a continuous piecewise quadratic finite element space on the fixed background mesh $\mesh$  denoted by $\mathcal{W}_h$. Then, to define our discrete domains $\Omega_{i,h}$ and $\Gamma_h$, we use a two-grid solution proposed by \cite{Gross2006, Gross2011}. This two-grid solution can be outlined as follows. Firstly, we project the piecewise quadratic level set function into a piecewise linear finite element space defined by
a hierarchically refined mesh $\mathcal{T}_{h/2}$ (\textit{i.e.} each edge is bisected). Let us call $\mathcal{I}(\phi_h)$ the projection of the piecewise quadratic level set defined by
\begin{equation}
\mathcal{I}(\phi_h(\mathbf{v})) = \phi_h(\mathbf{v}) \mbox{ for all nodes $\mathbf{v}$ in } \mathcal{T}_{h/2}.
\end{equation}
We then use this piecewise linear interpolation to determine the intersection between $\mathcal{I}(\phi_h)$ and the refined grid to obtain the piecewise linear approximation of $\Omega_{i,h}$ and $\Gamma_h$ as illustrated in Figure~\ref{fig: refined mesh approximation}.  \\
\paragraph{Interface normals} 
The normal pointing from $\Omega_1$ to $\Omega_2$ can be obtained using the level set function as 
\begin{equation}
n_{\Gamma}(x,t) = \frac{\nabla \phi(x,t)}{\norm{\nabla \phi(x,t)}}. 
\end{equation}
In the finite element approximation, the normal $n_{\Gamma}$ is obtained from the piecewise quadratic level set function through a $L^2$-projection onto the continuous piecewise linear space 
\begin{equation}
\mathcal{X}_h^d := \left\{ v_h \in [C^0(\Omega)]^d : \restr{v_h}{K} \in \mathcal{P}_1(K)\, \forall K \in \mesh \right\}.
\end{equation} 
Here, $d=2,3$ is the geometrical dimension. We determine the normal by finding $n_{\Gamma} \in \mathcal{X}_h^d $ such that for all $\delta n_{\Gamma}  \in \mathcal{X}_h^d$
\begin{equation}
(n_{\Gamma},\delta  n_{\Gamma})_{\Omega} = \left(\frac{\nabla \phi}{|\nabla \phi|},\delta  n_{\Gamma}\right)_{\Omega}. 
\label{equ: l2 projection normal}
\end{equation} 

\paragraph{Reinitialisation} To ensure that the level set function is close to the signed distance function, we reinitialise the level set function in regular time intervals using a fast marching scheme as described in~\cite{Gross2011,Gross2006}.

\subsubsection{Stabilised Cut Two-Phase Navier-Stokes Discretisation}

Let $ \mcV_h $ denote the space of piecewise linear polynomials
on $\mcT_h$, i.e. 
\begin{align}
  \mathcal{V}_h &= 
  \left\{ 
    v_h \in C^0(\Omega): \restr{v_h}{K} \in
\mathcal{P}_1(K)\,,\, \forall K \in \mcT_h
\right\}
\end{align}
and let 
\begin{equation}
\mathcal{\bfV}_h^n = \mathcal{V}_{h,1}^n \oplus \mathcal{V}_{h,2}^n \quad \mbox{with } \mathcal{V}_{h,i}^n = \left. \mathcal{V}_h \right|_{\mathcal{T}_{h,i}(t_n)}, i = 1,2.
\end{equation}
Choosing $p_h(t_n) \in \mathcal{\bfV}_h^n$ now means that  $p_h(t_n)$ consists of a
pair $(p_{h,1}^n, p_{h,2}^n)$, with $p_{h,i}^n \in  \mathcal{V}_{h,i}^n$, and
$p_h(t_n) \in \mathcal{\bfV}_h$ is double valued for elements that are both in $\mathcal{T}_{h,1}(t_n)$ and $\mathcal{T}_{h,2}(t_n)$, i.e. all elements in the 
$\delta$ neighbourhood $G_h^n$ of the moving interface. This allows $p_h$ to be discontinuous at
the interface (see Figure~\ref{fig: domains}) similar to XFEM methods. Analogously, we define the vector-valued  linear finite element
space $[\mathcal{\bfV}_h^n]^d$  for the velocity. Note that the velocity is also double valued for elements in $G_h^n$ and together with the zero jump condition, this yields a velocity which can exhibit kinks at the interface. Additionally, the definition of $\vel_{h,i}$, $i=1,2$ on the extended fictitious domains $\Omega_i^*$ will be used to evaluate the velocity of the previous time-step, $\vel_{h,i}^{n-1}$, in the updated domains $\Omega_{h,i}^n$. In this contribution, the extension of the velocity and the pressure from the physical domain, $\Omega_{h,i}^n$,  to the fictitious domain, $\Omega^*_i(t_n)$ will be achieved through facet stabilisation terms defined hereafter. \\
The stabilised finite element formulation of the two-phase Navier-Stokes problem reads: For $t_n \in \{ t_1, \dots t_{n_t}\}$,
  find $U_h := (\vel_h^n,p_h^n) \in  [\mathcal{\bfV}_h^n]^d \times \mathcal{\bfV}_h^n $ 
  such that for all $V_h:=(\velt_h,q_h) \in  [\mathcal{\bfV}_h^n]^d \times \mathcal{\bfV}_h^n$
  \begin{equation}
    A(U_h,V_h) + J(U_h,V_h) = L(V_h),
\label{equ: discrete Navier-Stokes}
  \end{equation}
  where
  \begin{align}
  \label{eq:A_h-definition}
    A(U_h,V_h) =&  \rho \left( \frac{\vel_h^{n}}{\Delta t}, \velt_h\right)_{\Omega} + \mathcal{N}(\bfbeta_h,\vel_h^{n},\velt_h) + a(\vel_h^{n},\velt_h) + a_n(\vel_h^{n},\velt_h) 
\nonumber \\ \phantom{=}& + b(p_h^{n},\velt_h) -  b(q_h,\vel_h^{n}), \\
  L(V_h) =&  \rho \left( \frac{\bfbeta_h}{\Delta t}, \velt_h \right)_{\Omega} +f_{\Gamma}(\velt_h), \\ 
J(U_h,V_h) =& j_{\beta}(\vel_h,\velt_h) +j_{div}(\vel_h,
 \, \velt_h) + j_p(p_h,q_h) + g_u(\vel_h,\velt_h) , \\
f_{\Gamma}(\velt_h)  =& -(\gamma \kappa \mathbf{n}_{\Gamma}, \mean{\velt_h})_{\Gamma_h^n}
\label{equ: two-phase Navier discretised}
  \end{align}
with given initial condition $\vel(\bfx,0) = \vel_0$.
Here, the integrals are evaluated as 
\begin{equation}
\left( \vel_h, \velt_h \right)_{\Omega} = \sum_{i=1}^2 \left( \vel_{h,i}, \velt_{h,i} \right)_{\Omega_{h,i}^n}. 
\end{equation}
And the penalty terms to enforce the interface and boundary conditions are scaling with the mesh size $h$ as follows: $\lambda_{\partial \Omega}(h) = \frac{\hat{\lambda}_{\partial \Omega}}{h}$,  $\lambda_{\Gamma}(h) = \frac{\hat{\lambda}_{\Gamma}}{h}$.
The stabilisation terms in $J$ consist of jump penalty terms over facets of the active meshes. Inspired by  \cite{Winkelmann2007,BurmanFernandezHansbo2006,Burman2007}, we use the following stabilisation terms
\begin{align}
j_{\beta}(\vel_h,\velt_h) &=  \gamma_u
 \sum\limits_{i=1}^2 \sum\limits_{F \in \mathcal{F}_I^{i}} h^2  \rho_i \norm{\bfbeta_h \cdot \bfn_F}_{\infty,F} \left( \jump{\nabla \vel_h}_{\bfn_F} ,\jump{\nabla \velt_h}_{\bfn_F} \right)_F , \label{equ: j_b} \\
  j_{div}(\vel_h,\velt_h) &=  \gamma_{div}
\sum\limits_{i=1}^2 \sum\limits_{F \in \mathcal{F}_I^{i}} h^2  \rho_i \norm{\bfbeta_h}_{\infty,K} \left( \jump{\nabla \cdot \vel_h } ,\jump{\nabla \cdot \velt_h } \right)_F , \\
j_p(p_h,q_h) &= \gamma_p
  \sum\limits_{i=1}^2  \sum\limits_{F \in \mathcal{F}_I^{i}} \frac{h^3}{ \eta_i \max(Re_K^i, 1)} \left( \jump{\nabla p_h}_{\bfn_F} ,\jump{\nabla q_h}_{\bfn_F} \right)_F ,
\end{align}
\begin{equation}
Re_K^i = \frac{h \rho_i \norm{\bfbeta_h}_{\infty,K}}{ \eta_i}.
\end{equation}
Additionally, we use a ghost penalty regularisation for the velocity, which is required at low Reynolds numbers 
\begin{align}
g_u(\vel_h,\velt_h) &=  \gamma_{gu}
 \sum\limits_{i=1}^2 \sum\limits_{F \in \mathcal{F}_G^{i}}  \eta_i h \left( \jump{\nabla \vel_h}_{n_F}  ,\jump{\nabla \velt_h}_{n_F} \right)_F .
\label{equ: g_u}
\end{align}
Here, $\jump{\nabla x}_{\bfn_F}$ denotes the normal jump of the quantity $x$
over the facet, $F$, defined as $\jump{\nabla x}_{\bfn_F} = \left. \nabla x
\right|_{K_F^+} {\bfn_F}  - \left. \nabla x
\right|_{K_F^-} {\bfn_F}$, where ${\bfn_F}$ denotes a unit
normal to the facet $F$ with fixed but arbitrary orientation. The
stabilisation parameters $\gamma_p$, $\gamma_u$, $\gamma_{div}$, $\gamma_{gu}$
have to be chosen large enough.
The term $j_{\beta}$ stabilises the convective terms which is necessary for high Reynolds numbers, $j_{div}$ ensures additional control on the incompressibility condition necessary at high Reynolds numbers, $j_p$ ensures the inf-sup stability for an equal order approximation of velocity and pressure and prevents ill-conditioning. The term $j_p$ is required for low and high Reynolds number flows. In contrast, the terms $j_{\beta}$ and $j_{div}$ can be omitted for low Reynolds number flows, i.e. $\gamma_u=\gamma_{div} =0$ can be chosen. The term $g_u$ ensures stability independent of the cut location and prevents ill-conditioning for the velocity.  \\ 
In addition to stabilising the numerical scheme, facet terms \eqref{equ: j_b} -- \eqref{equ: g_u} extend the velocity and the pressure from the physical domains $\Omega_{h,i}$ to the fictitious domains $\Omega^*_i(t_n)$. This enables us to define the velocity of the previous time step in the current domains as 
\begin{equation}
\bfbeta_h(\bfx) = 
\begin{cases}
\vel_{h,1}^{n-1}(\bfx) &\mbox{ for } \phi(\bfx,t_n) \leq 0, \\
\vel_{h,2}^{n-1}(\bfx) &\mbox{ for } \phi(\bfx,t_n) >0. 
\end{cases}
\end{equation}
To ensure that the domains $\Omega^{n+1}_{h,i}$ are contained in the fictitious domains of $\Omega^*_i(t_n)$, we choose the $\delta$ neighbourhood as the domain formed by all intersected elements and their neighbouring elements. That means a band of one cell layer on each side of the intersected elements resulting in a $\delta$-neighbourhood of three cell bands (see Figure~\ref{fig: extended velocity}). The interface strip thickness can also be chosen dependent on the velocity at the interface, e.g. $\delta \geq \norm{\vel \cdot \bfn_{\Gamma}}_{\infty,\Gamma_h} \Delta t$ \cite{Lehrenfeld2018} to ensure that  $\Omega^{n+1}_{h,i} \subset \Omega^*_i(t_n)$ and the level set function can be used to define the $\delta$-neighbourhood (if the level set function is sufficiently close to a signed distance function). However, as we are facing time step restrictions due to the decoupling of the Navier-Stokes equation and the level set advection as well as from the linearisation of the convective term, the 3-cell band $\delta$-neighbourhood is numerically shown to be thick enough.   

\subsubsection{Discretised Level Set Advection}
As mentioned above we discretise the level set function using a continuous quadratic finite element space, $\mathcal{W}_h$. We stabilise the advection problem~\eqref{equ: level set time discretisation} in space with streamline diffusion (SUPG). The discretised advection problem reads: Find $\phi_h^{n} \in \mathcal{W}_h$, such that for all $\delta \phi \in \mathcal{W}_h$
\begin{equation}
a_{\phi}(\phi_h^{n},\delta \phi) = l_{\phi}(\delta \phi)
\label{equ: level set advection}
\end{equation}
with 
\begin{equation}
\begin{aligned}
a_{\phi}(\phi_h^{n},\delta \phi) &= \left( \frac{\phi_h^{n}}{\Delta t}+ \theta \vel_h^{n} \cdot \nabla \phi_h^{n}, \delta \phi + \tau_{SD} (\vel_h^n \cdot \nabla \delta \phi) \right)_{\Omega}, \\
l_{\phi}(\delta \phi) &=  \left( \frac{\phi_h^{n-1}}{\Delta t}+ (\theta-1) \bfbeta_h \cdot \nabla \phi_h^{n-1}, \delta \phi + \tau_{SD} (\vel_h^n \cdot \nabla \delta \phi) \right)_{\Omega}
\label{equ: level set SUPG}
\end{aligned}
\end{equation}
with the streamline diffusion parameter \cite{Hansbo2016}
\begin{equation}
\tau_{SD} = 2 \left( \frac{1}{\Delta t^2} + \frac{| \vel_h^n \cdot \vel_h^n | }{h^2} \right)^{-\frac{1}{2}}
\end{equation}
and initial condition $\phi_h^0 =\phi_0$ and Dirichlet conditions on inflow boundaries. 

\subsection{Stabilised Scheme for the Surface Tension Force}
In this section, we introduce a scheme to smoothen the curvature approximation of the fluid interface. This smoothing is used to reduce the magnitude of spurious velocity modes which arise from errors in the curvature computation. This is particularly important for flows in which surface tension forces dominate viscous forces as these spurious velocities can be the source of numerical instabilities. \\
The surface tension force term 
\begin{equation}
f_{\Gamma}(\velt_h) = -(\gamma \kappa \mathbf{n}_{\Gamma}, \mean{\velt_h})_{\Gamma_h^n}
\label{equ: surface tension form}
\end{equation}
contains the mean curvature $H:=\kappa \bfn_{\Gamma}$ which can be reformulated in terms of the Laplace Beltrami operator as
(see e.g.  \cite{Hansbo2015},\cite{Hysing2006}, \cite{Gross2006})
\begin{equation}
H= \kappa \bfn_{\Gamma} = - \Delta_{\Gamma} id_{\Gamma}.
\end{equation}
Here, $id_{\Gamma}$ denotes the identity mapping on $\Gamma$ and $\Delta_{\Gamma} = \nabla_{\Gamma} \cdot \nabla_{\Gamma}$ is the Laplace-Beltrami operator with $\nabla_{\Gamma}$ denoting the tangential gradient given by 
\begin{equation}
\nabla_{\Gamma} = P_{\Gamma} \cdot \nabla, \quad
P_N = I - \bfn_{\Gamma} \otimes \bfn_{\Gamma},
\end{equation}
where $I$ is the identity matrix and $\otimes$ is the outer product. 
Note that, $ \nabla_{\Gamma}id_{\Gamma} = P_{\Gamma} \cdot \nabla id_{\Gamma} = P_{\Gamma} I = P_{\Gamma}$.\\
 \\
We will now replace the mean curvature vector $H = \kappa \bfn_{\Gamma}$ with a smoothened discrete curvature \cite{Hansbo2016, Hansbo2015,Cenanovic2017}, $H_h$,  as detailed in the following. Let $[\mathcal{H}^n_h]^d$ denote the vector-valued continuous piecewise linear finite element space over  the domain formed by the band of all intersected elements at time $t=t_n$. Then, we determine a smoothened mean curvature vector $H_h \in [\mathcal{H}^n_h]^d$ by solving
\begin{equation}
\begin{aligned}
\left( H_h, \delta H \right) + J_H(H_h,\delta H) &=  ( P_{\Gamma}, \nabla_{\Gamma} \delta H )_{\Gamma_h^n}, \quad \forall  \delta H_h \in [\mathcal{H}^n_h]^d, \\
J_H(H_h,\delta H) &= \sum_{F \in \mathcal{F}_{int}(t_n)} \gamma_H \left( \jump{\nabla H_h}_{\bfn_F}, \jump{\nabla \delta H}_{\bfn_F} \right)_F.
\end{aligned}
\label{equ: mean curvature stabilised}
\end{equation}
Here, $\gamma_H$ is a positive penalty parameter. This stabilised mean curvature vector $H_h$ enters the right hand side of the two-phase Navier-Stokes problem \eqref{equ: two-phase Navier discretised} in the surface tension force term as  
\begin{equation}
f_{\Gamma}(\velt_h) = - (\gamma H_h, \mean{\velt_h})_{\Gamma_h^n}. 
\label{equ: stabilised curvature force term}
\end{equation}


\section{Coupled Formulation}
Algorithm~\ref{alg: cutfem algorithm} summarises the two-phase flow cut finite element scheme presented in the previous sections. 
\begin{algorithm}
    \caption{CutFEM Two-Phase Flow Algorithm}
    \label{alg: cutfem algorithm}
    \begin{algorithmic}[1] 
    \State Set $t=0$, $\vel_h^0=\vel_0$, $\phi_h^0 = \phi_0$.
    \While{$t\leq T$} \Comment{$T$ is the final time}
   \State Determine $\Omega_h^n$ and $\Gamma_h^n$ through intersection computations of zero level set with background mesh.   
\State Compute normal $\bfn_{\Gamma}$ using \eqref{equ: l2 projection normal}.
\State Compute the stabilised mean curvature vector $H_h$ using \eqref{equ: mean curvature stabilised}.
  \Procedure{two-phase Navier-Stokes}{$\vel_h^{n-1}$, $\bfn_{\Gamma}$, $H_h$, $\Omega_h^n$, $\Gamma_h^n$} 
	\State Solve the Navier-Stokes problem \eqref{equ: discrete Navier-Stokes}.
	 \State \textbf{return}	$\vel^{n}_h$, $p_h^n$. 
  \EndProcedure
  \Procedure{Level Set Advection}{$\vel_h^n$}
\State Solve advection problem~\eqref{equ: level set SUPG} for level set. 
\State Reinitialise level set function using fast marching scheme \cite{Gross2006,Gross2011}.
\State Apply mass correction using volume correction \cite{Gross2006,Gross2011}.
\State \textbf{return}	$\phi_{h}^{n}$. 
\EndProcedure
       \State $t =t+\Delta t$, $\phi_h^{n-1} = \phi_{h}^{n}$, $\vel^{n-1}_h = \vel^{n}_h$.
    \EndWhile\label{while}
    \end{algorithmic}
\end{algorithm}

\section{Numerical Results}

The proposed two-phase flow solver is implemented in the CutFEM library \cite{BurmanClausHansboEtAl2014} based on FEniCS \cite{fenics2015,Hale2017}. 

\subsection{Rising Bubble Benchmark}
\begin{figure}
\centering
\includegraphics[width=.3\textwidth]{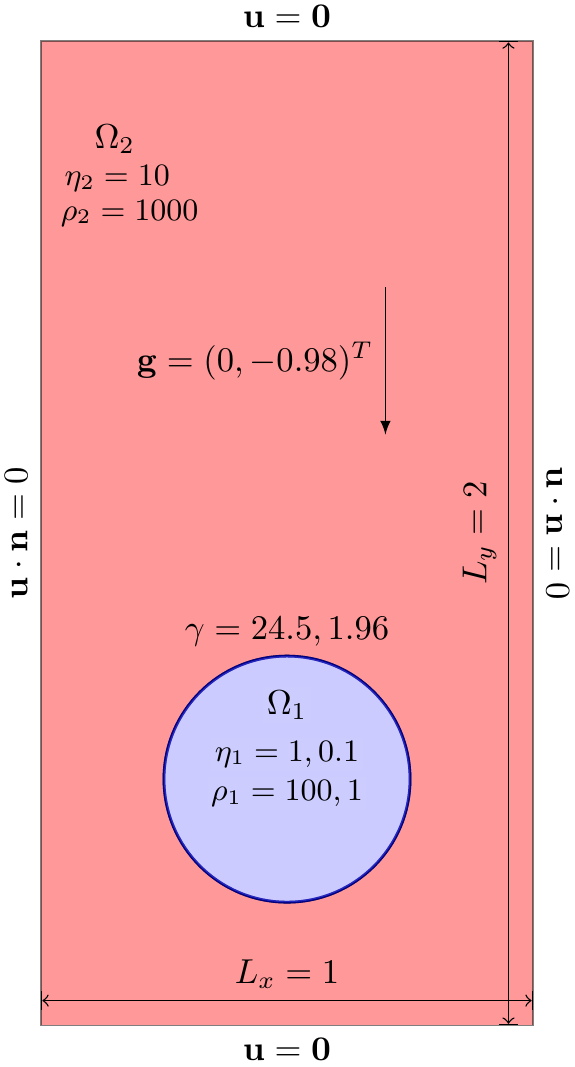}
\caption{Schematics of rising bubble benchmark.}
\label{fig: schematics rising bubble benchmark}
\end{figure}

\begin{figure}
\centering
\subfloat[Bubble shape at time $t=3$.]{\includegraphics[width=.5\textwidth]{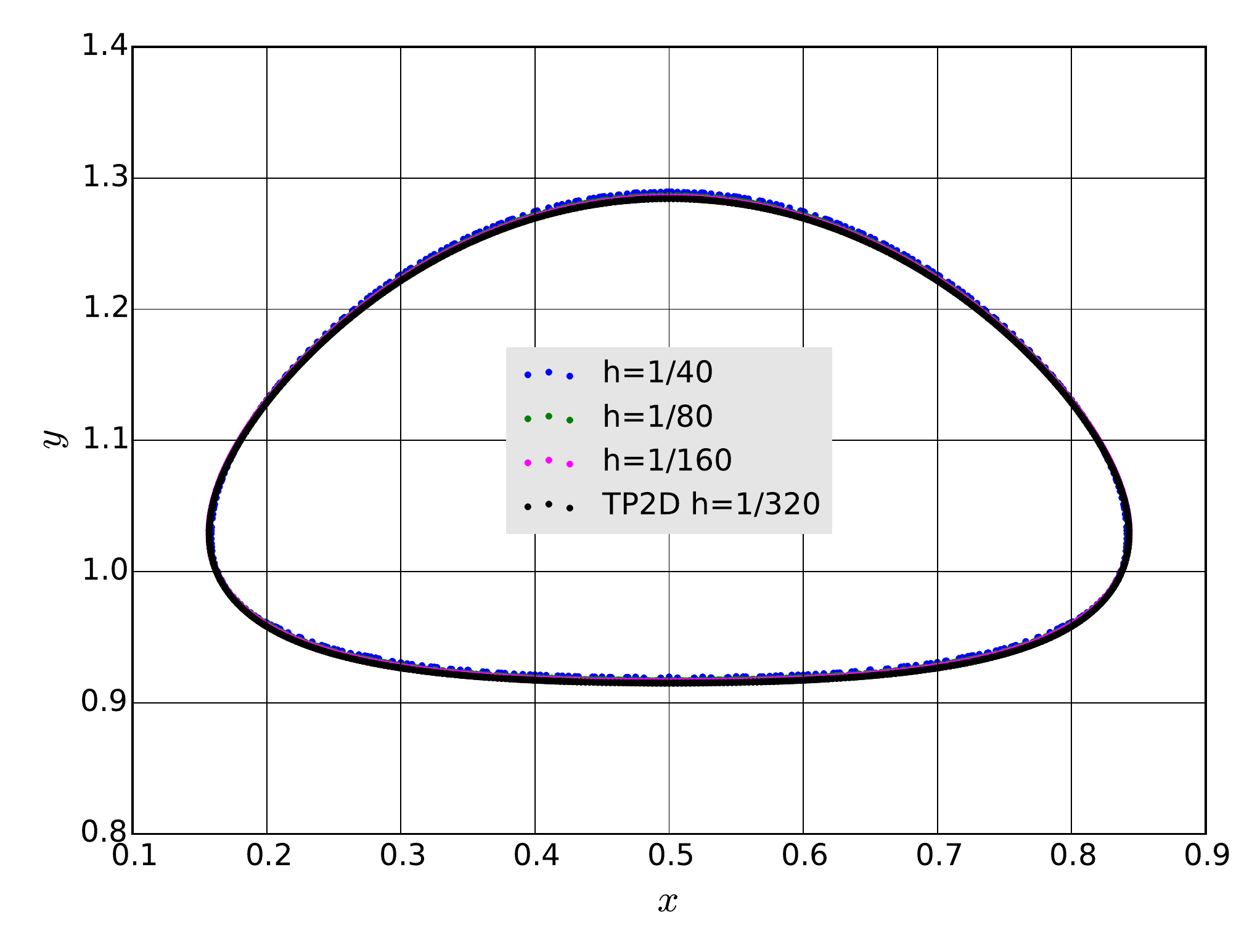}}
\subfloat[Center of Mass.]{\includegraphics[width=.5\textwidth]{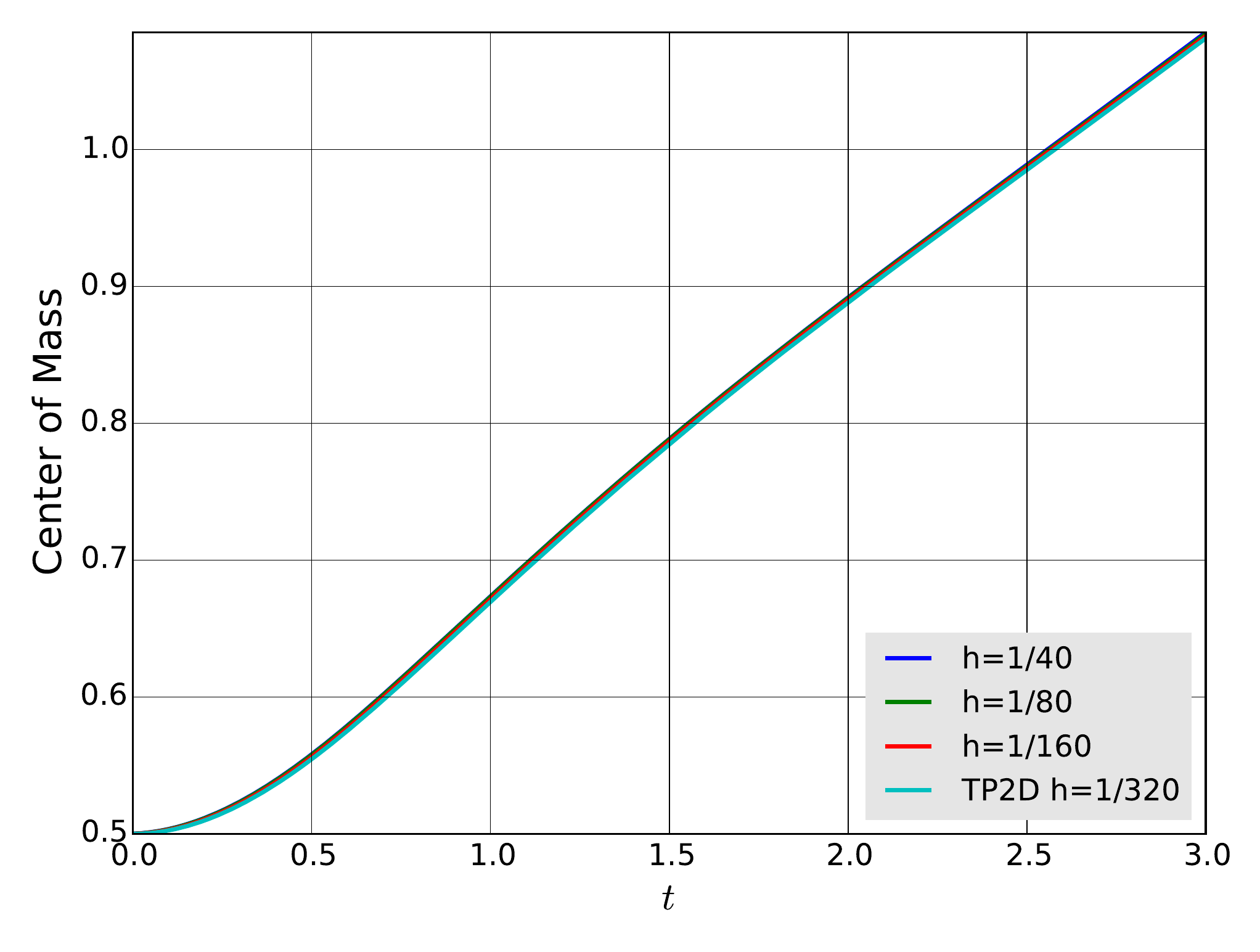}}\\
\subfloat[Rise Velocity.]{\includegraphics[width=.5\textwidth]{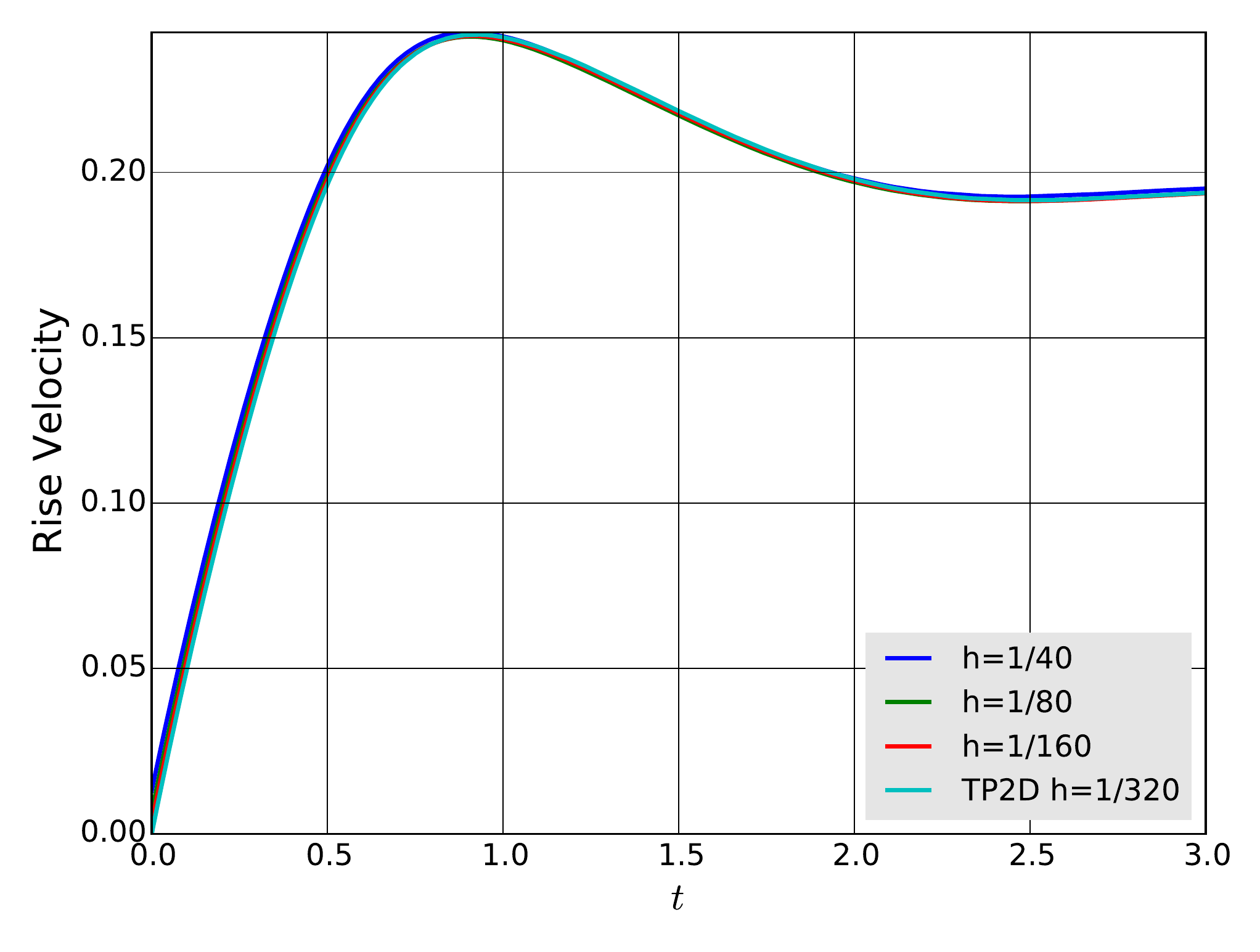}}
\subfloat[Circularity.]{\includegraphics[width=.5\textwidth]{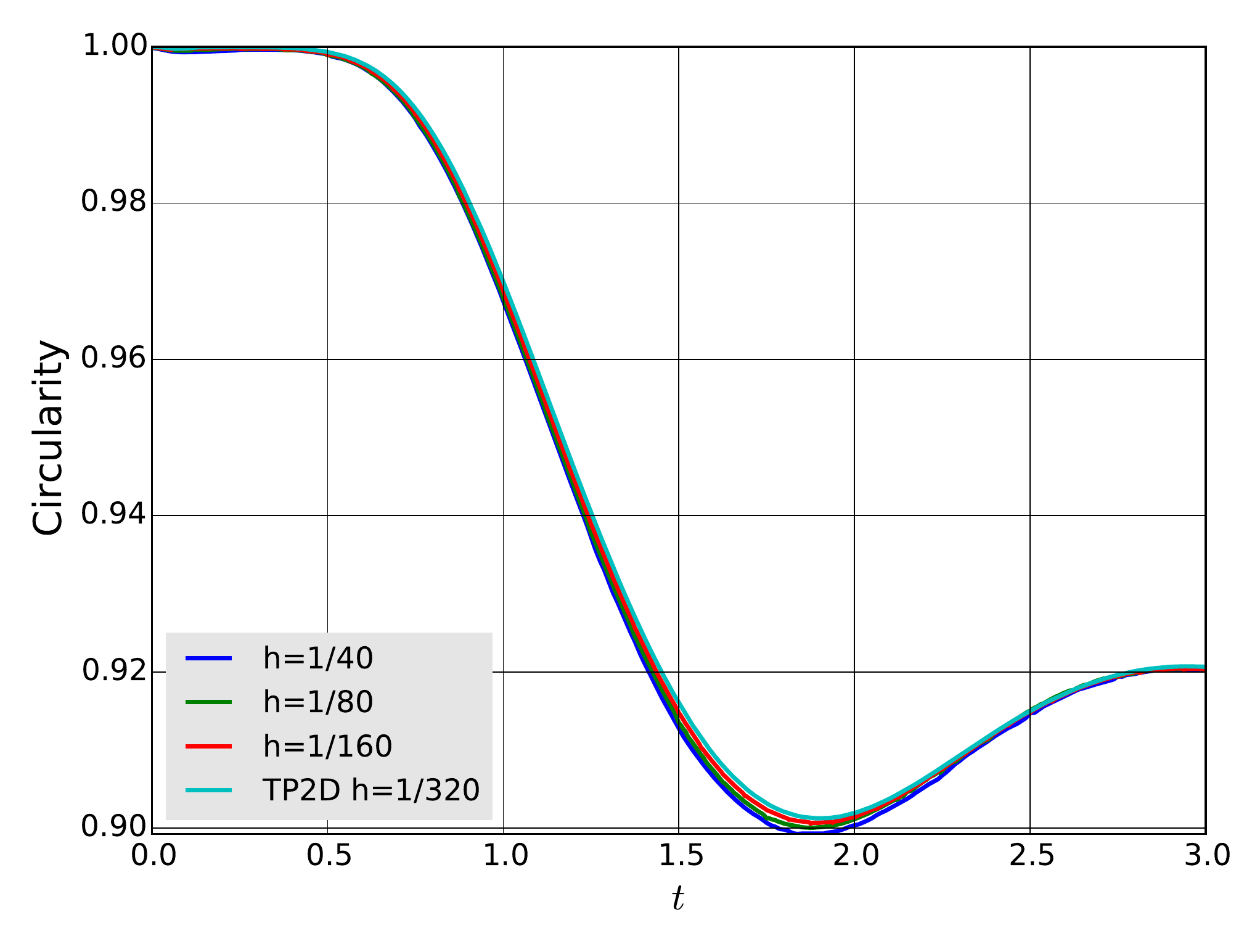}}\\
\caption{Numerical results of rising bubble benchmark for test case 1.}
\label{fig: rising bubble results case 1}
\end{figure}

\begin{figure}
\centering
\subfloat[Bubble Shape with curvature stabilisation at $t=1.9$.]{\includegraphics[width=.6\textwidth]{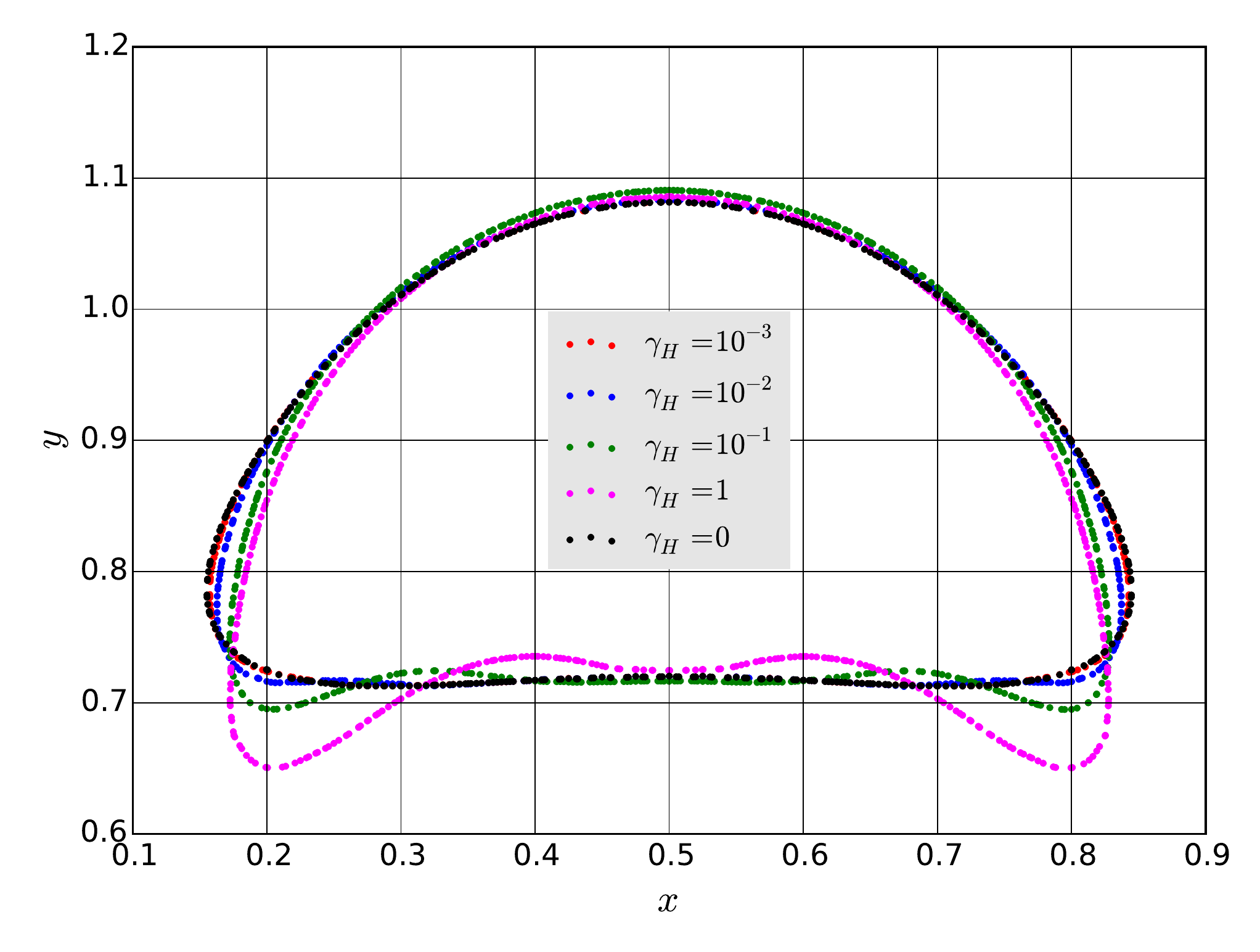}} \\
\subfloat[Velocity magnitude for Laplace Beltrami operator.]{\includegraphics[width=.39\textwidth]{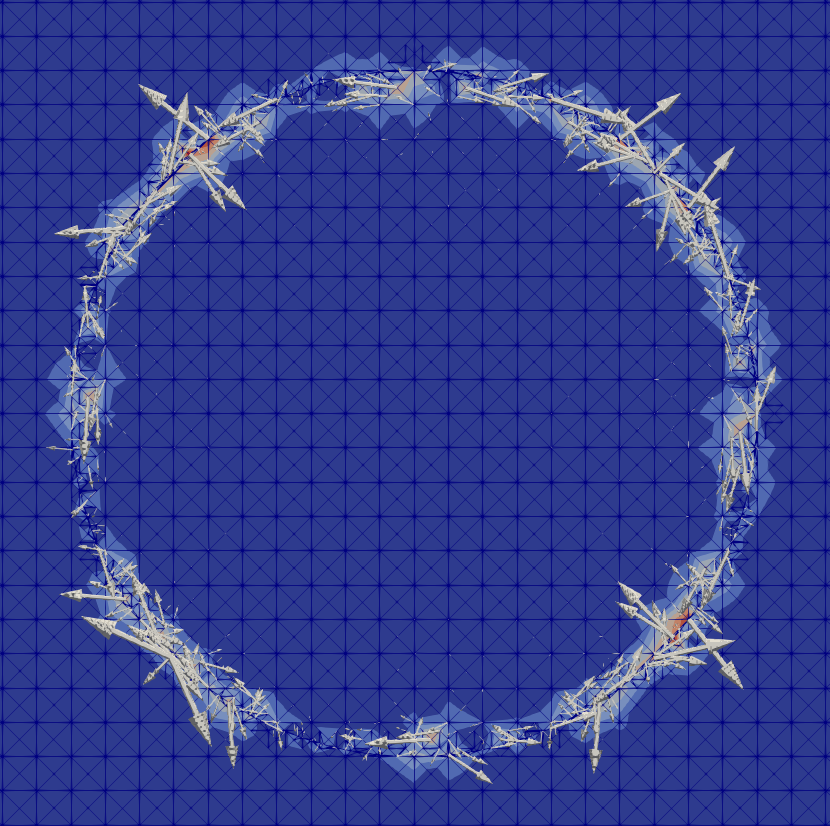}
\includegraphics[width=.09\textwidth]{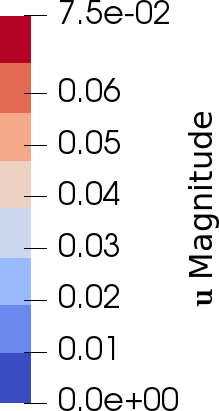}} 
\hspace{.01\textwidth}
\subfloat[Velocity magnitude for stabilised mean curvature.]{\includegraphics[width=.39\textwidth]{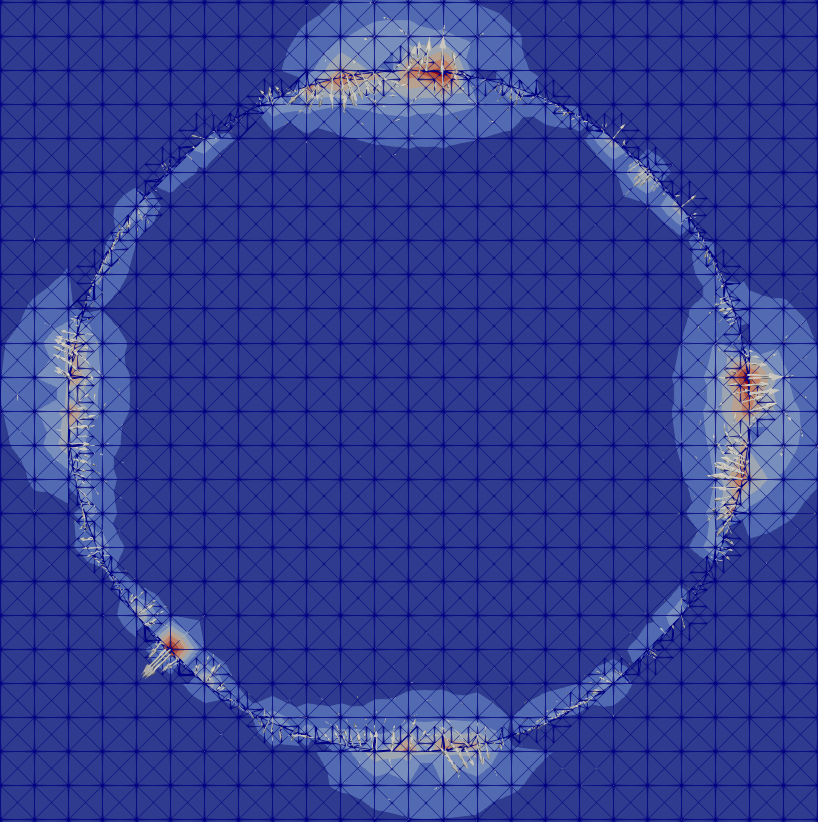}
\includegraphics[width=.09\textwidth]{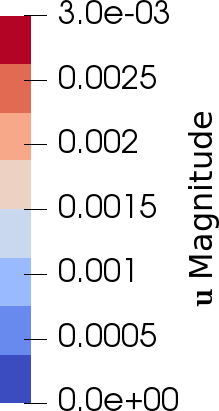}} 
\caption{Effect of curvature stabilisation on the bubble shape (a) and the reduction in spurious velocities for rising bubble test case 1 with zero gravity at time $t=0.125$ from Laplace Beltrami operator ($\gamma_H=0$) to stabilised curvature computation ($\gamma_H=0.001$) (b),(c).}
\label{fig: curvature stabilisation}
\end{figure}

\begin{figure}
\subfloat[]{
\includegraphics[width=.44\textwidth]{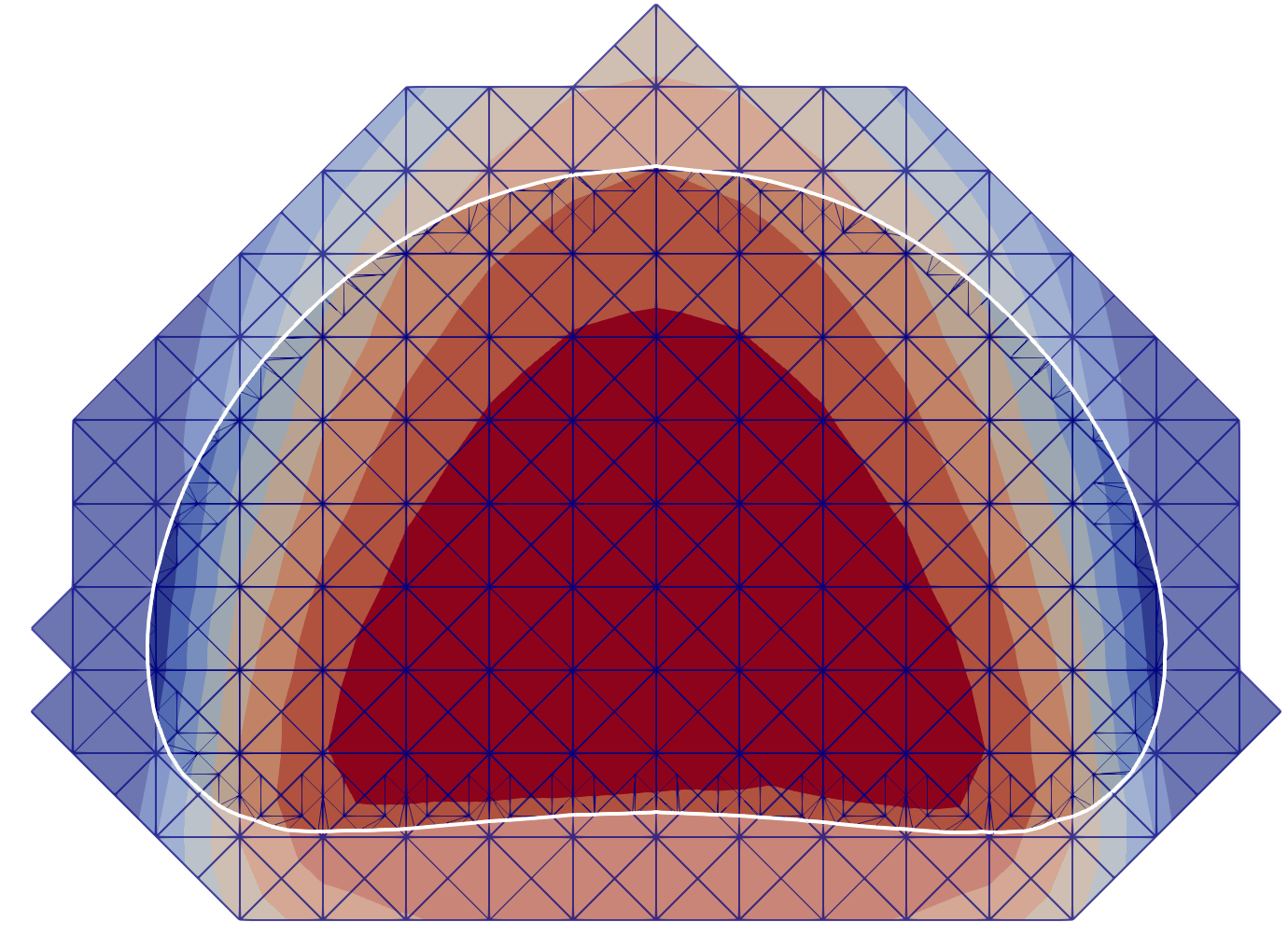}
\includegraphics[width=.1\textwidth]{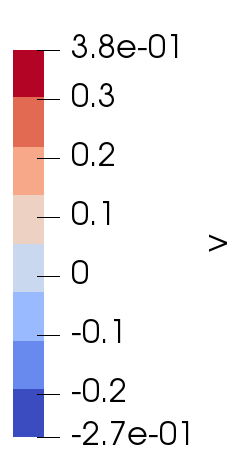}}
\subfloat[]{
\includegraphics[width=.46\textwidth]{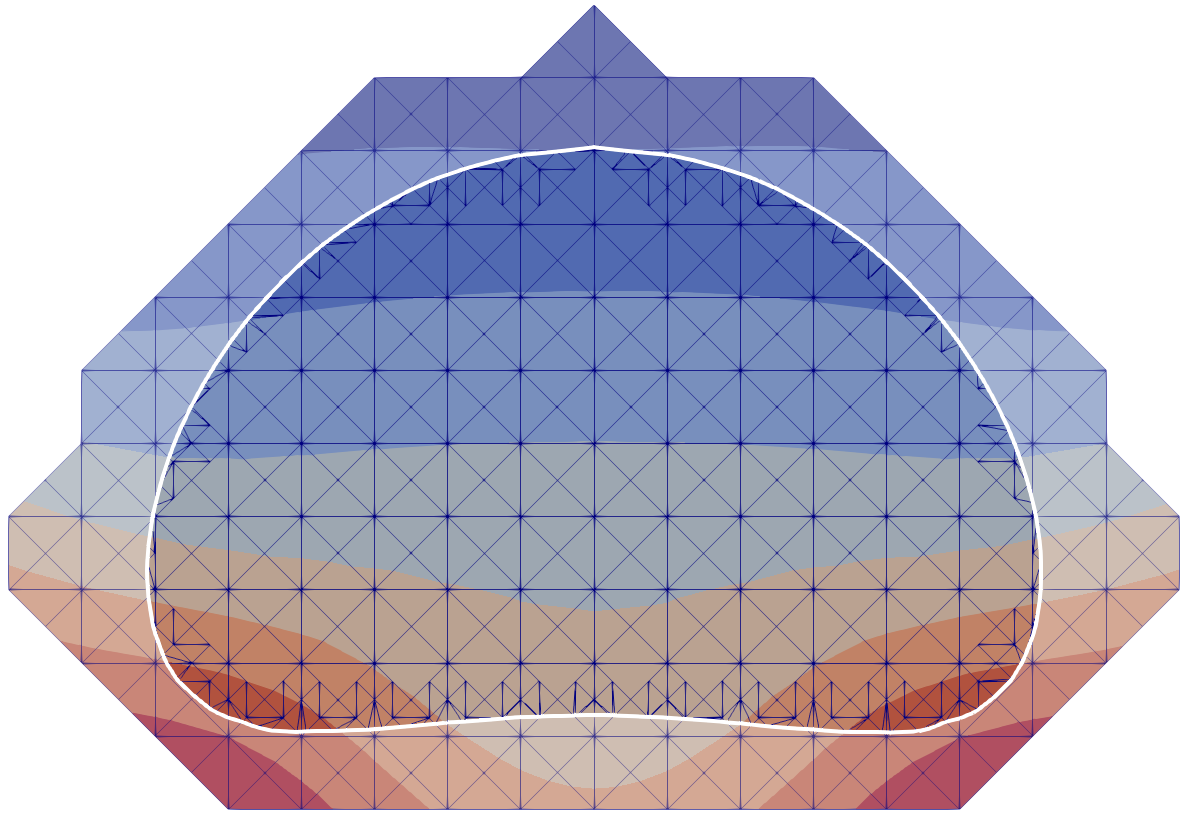}
\includegraphics[width=.1\textwidth]{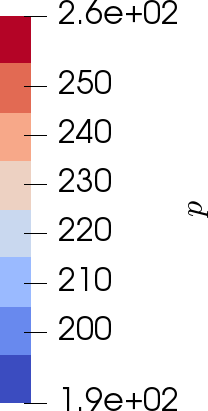}}

\caption{Extended $y$-component of the velocity (a) and the pressure (b) for test case 1 for a coarse mesh with $h=1/20$ at time $t=1.225$.}
\label{fig: extended velocity}
\end{figure}

\begin{figure}
\subfloat[]{
\includegraphics[width=.5\textwidth]{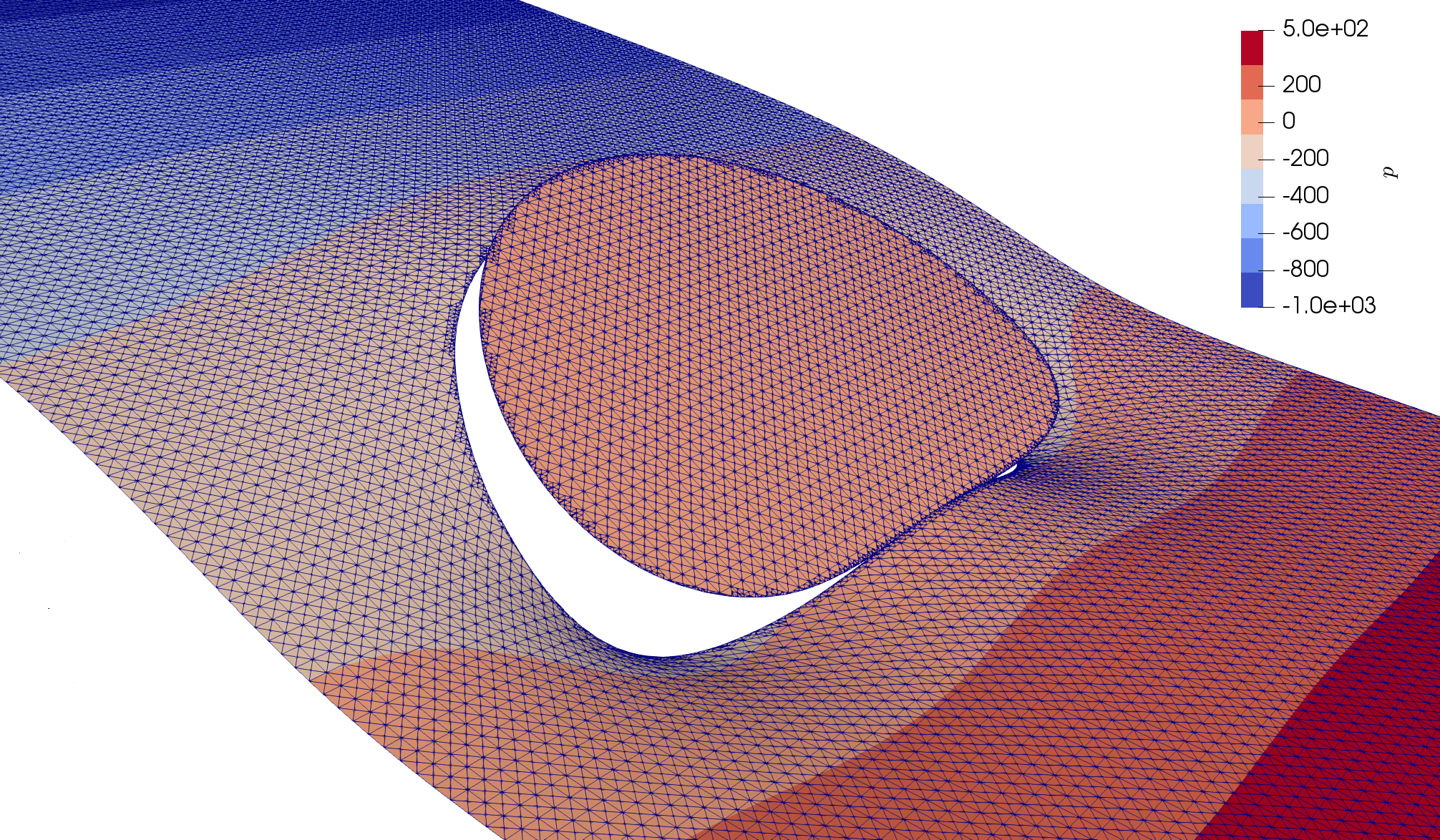}}
\subfloat[]{
\includegraphics[width=.5\textwidth]{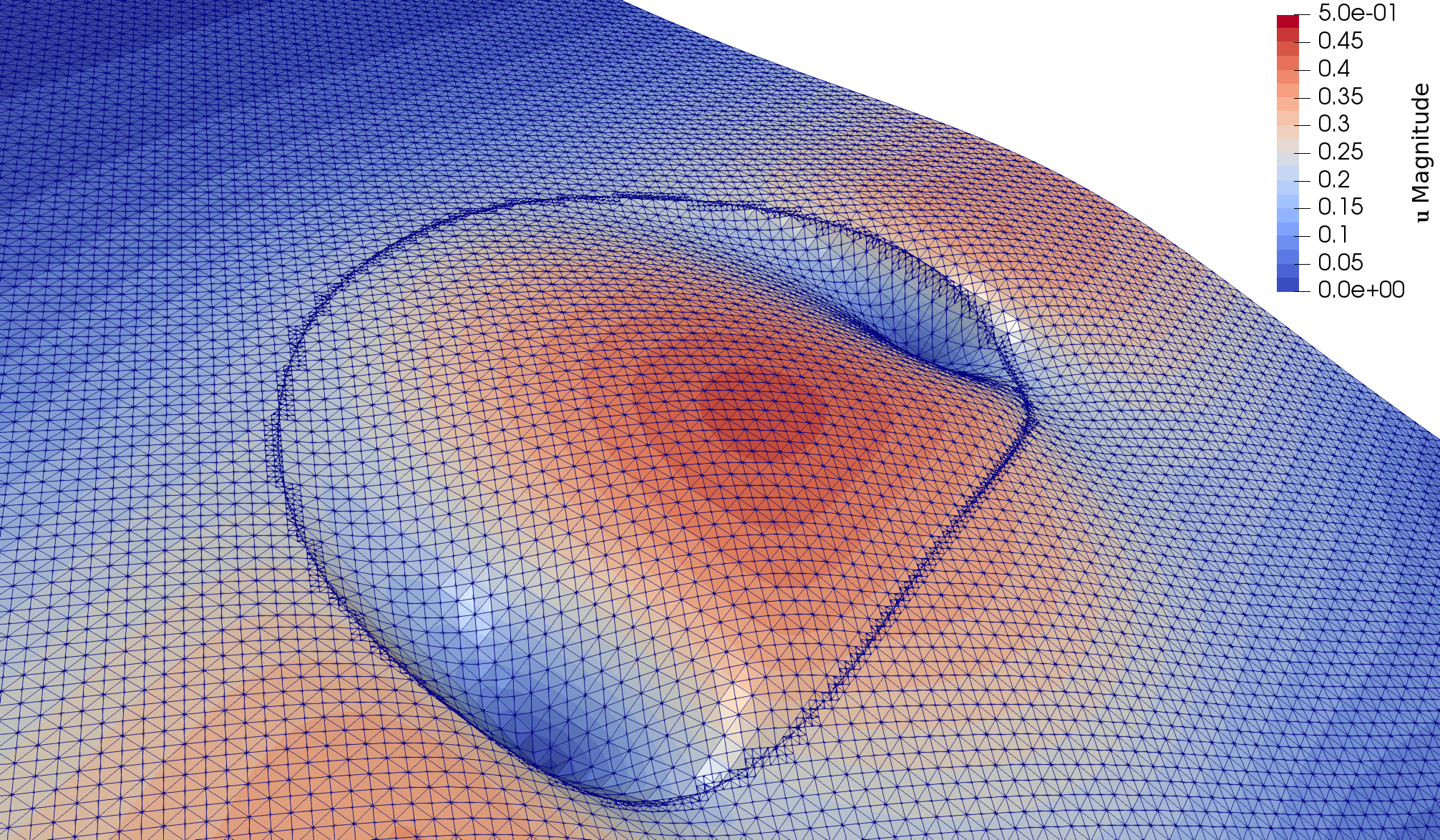}}
\caption{The cut finite element scheme is capable of representing kinks in the velocity and jumps in the pressure at the interface. Here, the pressure jump (a) and the velocity kink in the magnitude of the velocity (b)  is shown for test case 1 at time $t=1$ for $h=1/80$.}
\label{fig: pressure jump}
\end{figure}

\begin{figure}
\centering
\subfloat[Center of Mass.]{\includegraphics[width=.5\textwidth]{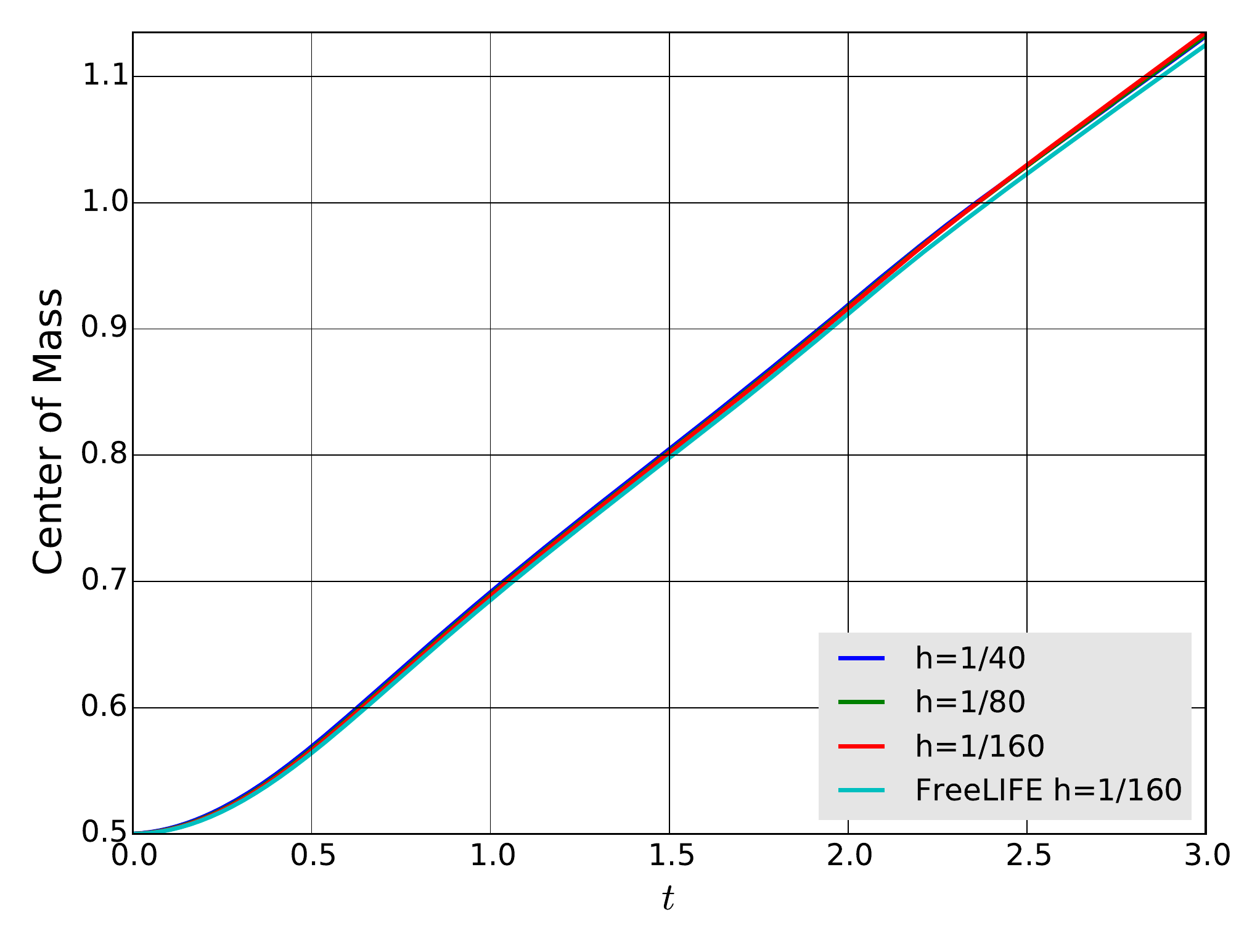}}
\subfloat[Rise Velocity.]{\includegraphics[width=.5\textwidth]{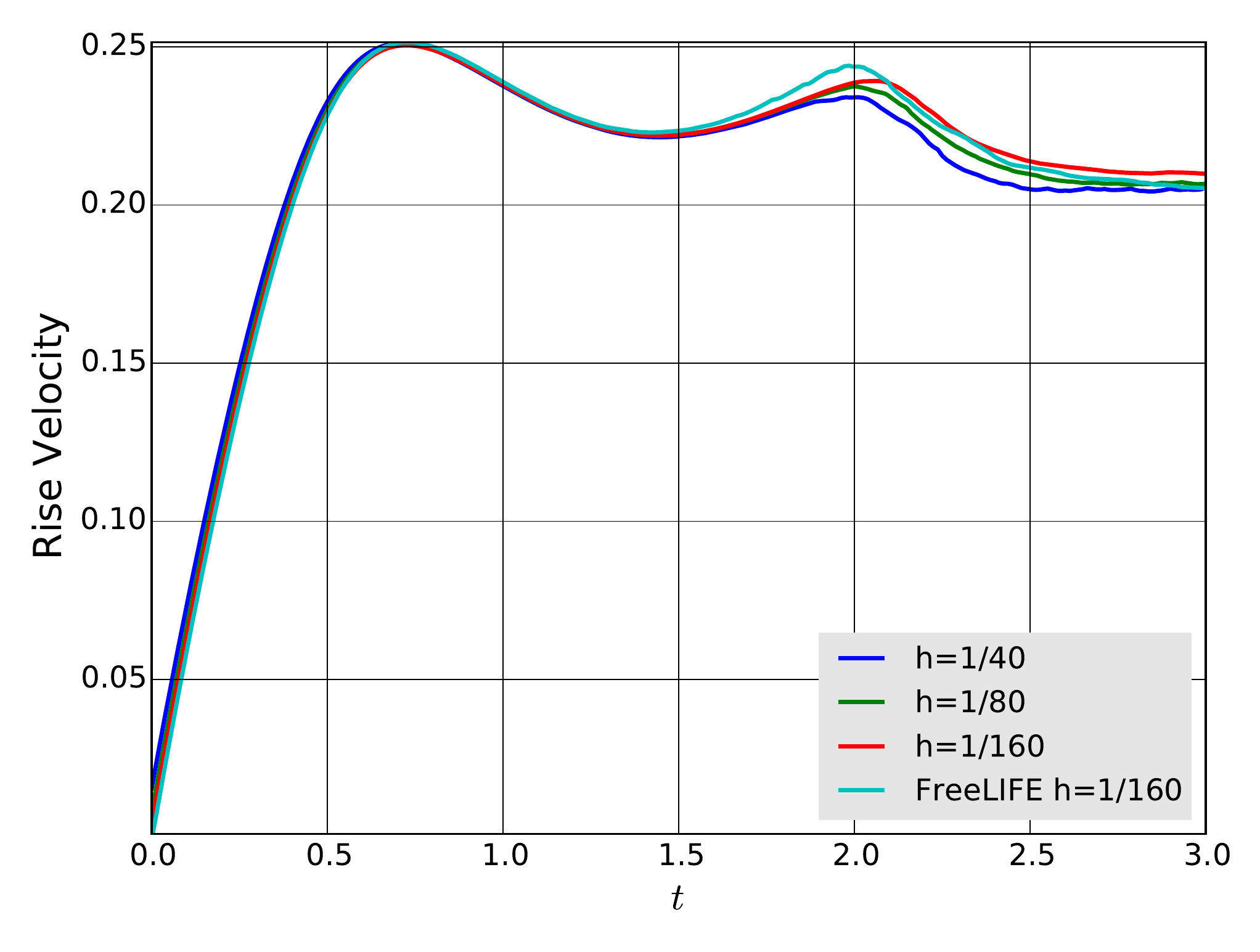}}\\
\subfloat[Circularity.]{\includegraphics[width=.5\textwidth]{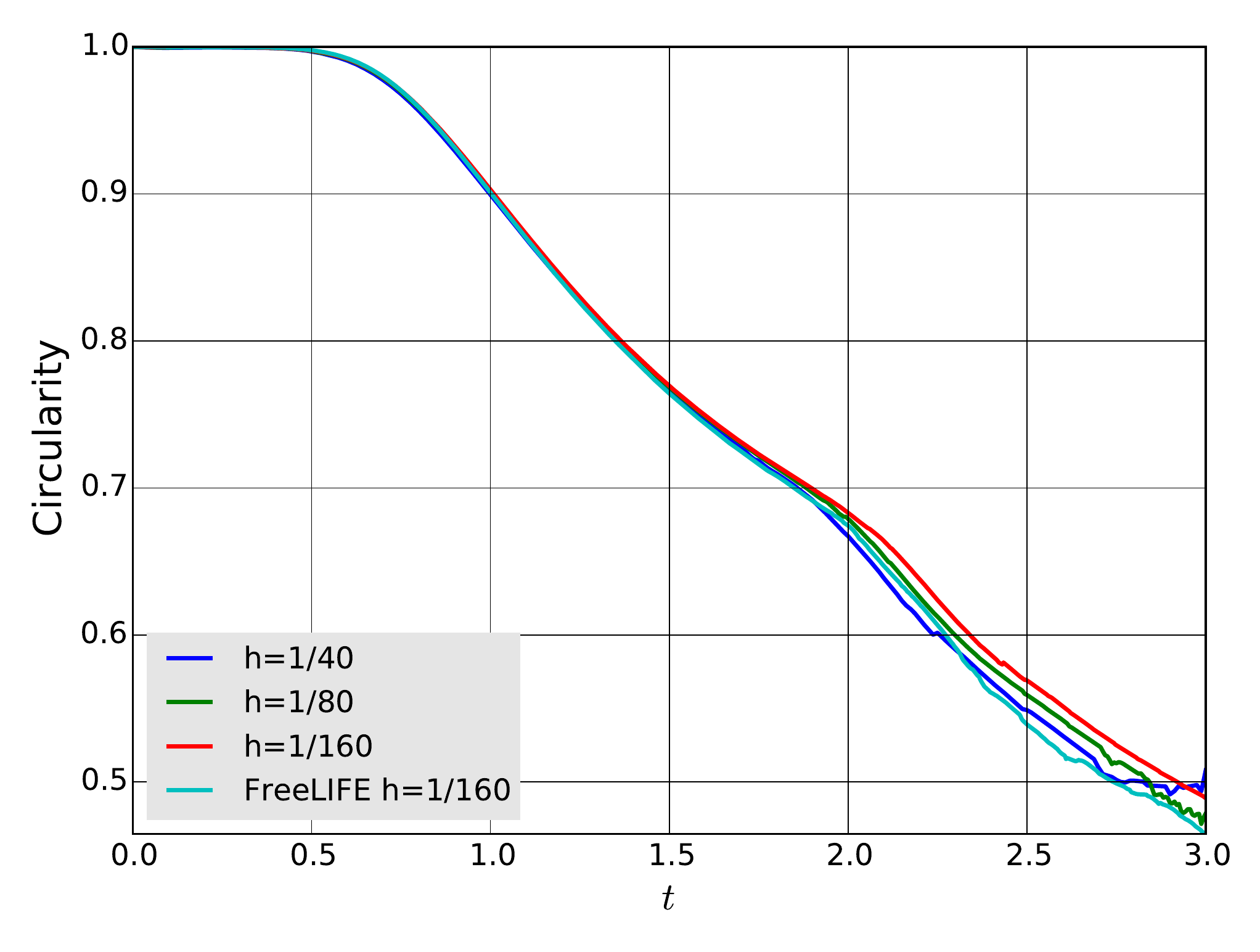}}
\caption{Center of mass, rise velocity and circularity for test case 2 for mesh refinement in comparison with FreeLIFE $h=1/160$.}
\label{fig: rising bubble results case 2}
\end{figure}

\begin{figure}
\centering
\subfloat[Bubble shape with mesh refinement.]{\includegraphics[width=.5\textwidth]{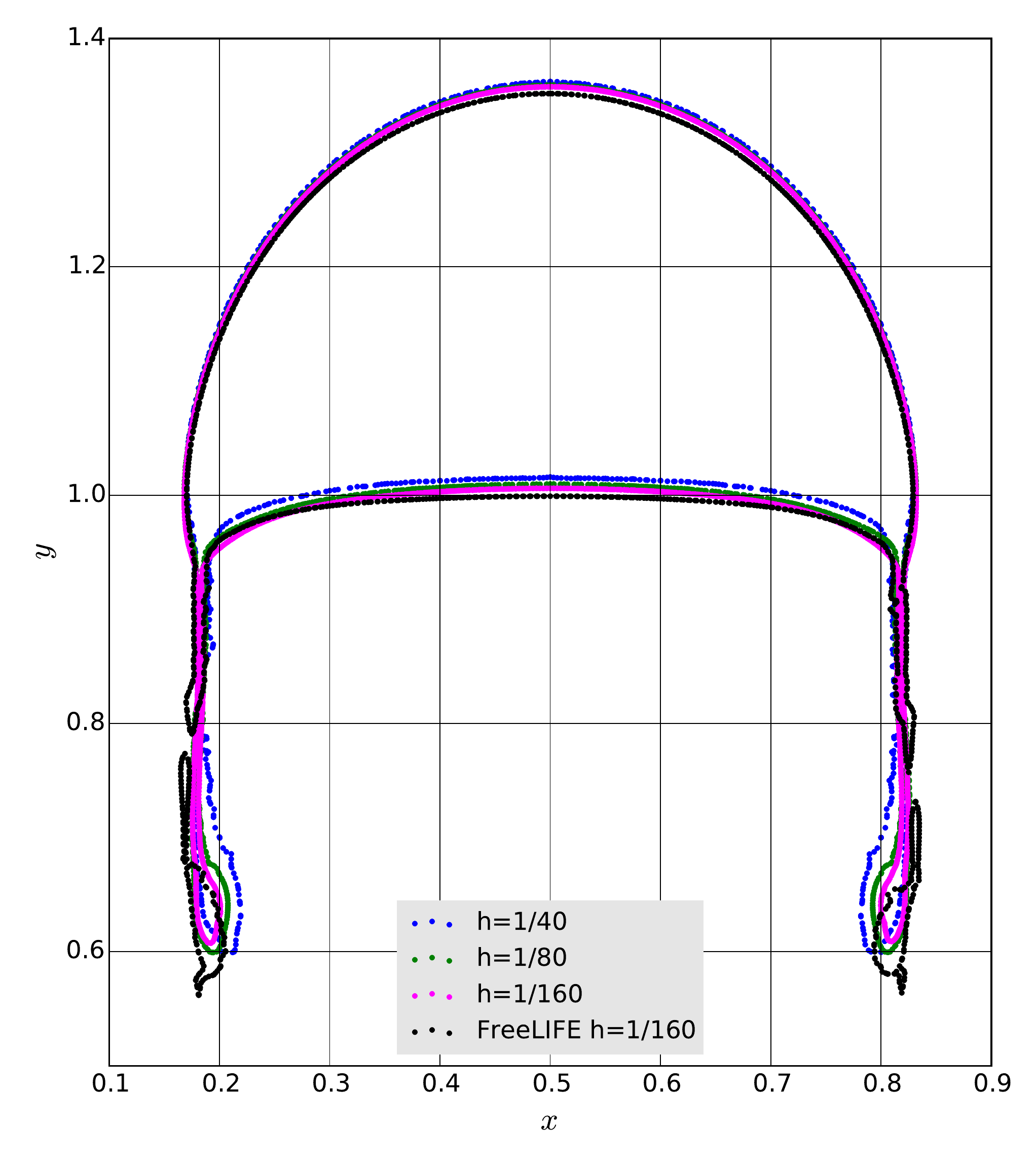}} 
\subfloat[Bubble shape with varying $\gamma_u$.]{\includegraphics[width=.5\textwidth]{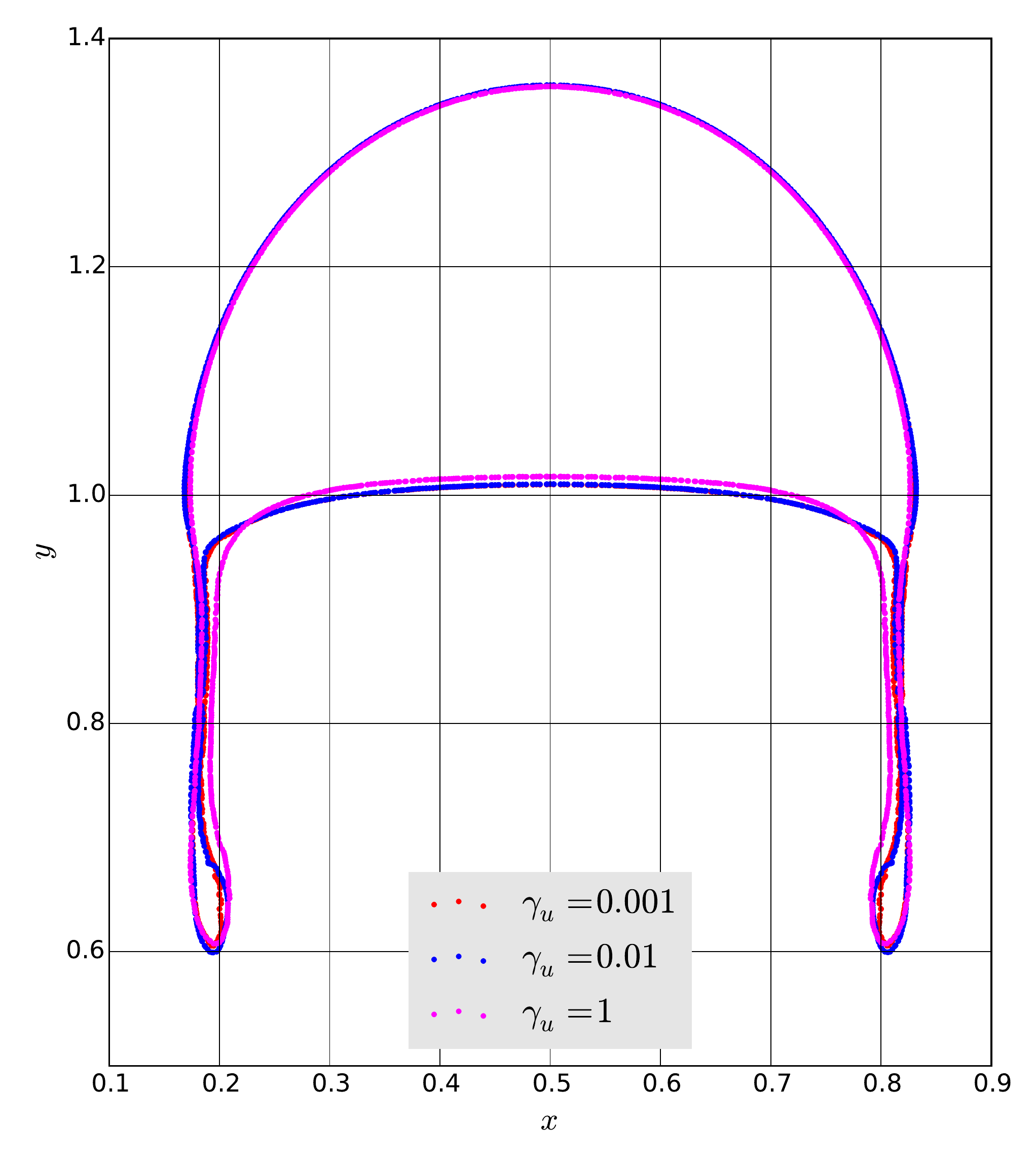}}
\caption{Bubble shapes for test case 2 with mesh refinement and with varying velocity stabilisation parameter $\gamma_u=\gamma_{div}=\gamma_{gu}$ for mesh size $h=1/80$ at $t=3.0$.}
\label{fig: rising bubble results case 2 shape}
\end{figure}

In this section, we validate our numerical scheme on the well known two-dimensional rising bubble benchmark presented in \cite{Hysing2009}. In this benchmark, a circular gas bubble with radius $r=0.25$ and centre point $\bfx_m = (0.5,0.5)^T$ is surrounded by a higher density fluid in a container occupying domain $\Omega = [0,1] \times [0,2]$. Due to buoyancy forces caused by gravity, $\bfg = (0,-0.98)^T$, the lighter gas bubble rises in the container. At the top and bottom boundary, we set homogenous Dirichlet conditions for the velocity and at the sides, we set slip conditions, i.e. $\vel \cdot \bfn =0$ (see schematics in Figure~\ref{fig: schematics rising bubble benchmark}). We choose a timestep of $\Delta t = h_{min}/2$, where $h_{min}$ is the smallest mesh size $h$ in the background mesh $\mesh$ and compute the solution over time interval $[0,3]$. As in the original benchmark, we study two cases as summarised in Table~\ref{tab: rising bubble cases}. In the first test case the rising bubble deforms into an elliptical shape and in the second test case the bubble develops long liquid strands. Like in the original benchmark, we evaluate the results in terms of the mean rise velocity 
\begin{equation}
V_c = \frac{\int_{\Omega_1} u_y \, dx}{\int_{\Omega_1} 1 \, dx},
\end{equation}
where $u_y$ is the $y$-component of the velocity $\vel$; the $y$-component of the center of mass 
\begin{equation}
Y_c = \frac{\int_{\Omega_1} y \, dx}{\int_{\Omega_1} 1 \, dx}
\end{equation}
and the degree of circularity \cite{Wadell1933}
\begin{equation}
\bar{c} =
\frac{2 \sqrt{\pi \int_{\Omega_1} 1 \, dx}}{\int_{\Gamma} 1 \, ds},
\end{equation}
which compares the circumference of the bubble to the circumference of a circle with the same area as the bubble. \\ 
For both test cases, we set the penalty terms to $\gamma_u=0.01$, $\gamma_p=0.1$,  $\gamma_{div}=0.01$, $\gamma_{gu}=0.01$, $\gamma_H = 0.001$, $\hat{\lambda}_{\partial \Omega}=10.0$, $\hat{\lambda}_{\Gamma}=10.0$ unless specified otherwise. We investigate the results in comparison to the benchmark values for mesh sizes $h=\{1/40, 1/80, 1/160\}$. For test case 1, we compare our results to the finest mesh results ($h=1/320$)  of the code TP2D \cite{Turek1999} and for test case 2, we compare our results to the finest mesh results ($h=1/160$) of FreeLIFE \cite{Parolini2005} from the rising bubble benchmark \cite{Hysing2009}.  
\paragraph{Test case 1} Figure~\ref{fig: rising bubble results case 1} shows that we achieve excellent agreement between the reference solution of TP2D and our cut finite element Navier-Stokes scheme for the final bubble shape at $t=3$, for the evolution of the centre of mass, for the rise velocity and for the circularity. \\ 
Figure~\ref{fig: curvature stabilisation} shows the positive and negative effects of the smoothened curvature computation~\eqref{equ: mean curvature stabilised}. Considering the same material parameters as in test case 1 and setting the gravity to zero reveals well known spurious velocities caused by the surface tension force term~\eqref{equ: surface tension form}. Compared to the direct application of the Laplace Beltrami operator (see Figure~\ref{fig: curvature stabilisation}(b))
\begin{equation}
f_{\Gamma}(\velt_h) = -(\gamma \kappa \mathbf{n}_{\Gamma}, \mean{\velt_h})_{\Gamma_h^n} = -(\gamma \nabla_{\Gamma}id_{\Gamma}, \mean{\nabla_{\Gamma} \velt_h})_{\Gamma_h^n},
\label{equ: Laplace Beltrami force term}
\end{equation}
the smoothened mean curvature force term~\eqref{equ: stabilised curvature force term} helps control the magnitude of the spurious velocities (see Figure~\ref{fig: curvature stabilisation}(c)). However, if the curvature stabilisation parameter $\gamma_H$ is chosen too large the smoothened curvature approximation can change the shape of the bubble significantly. As shown in Figure~\ref{fig: curvature stabilisation}(a), the bubble is significantly more deformed for high values of $\gamma_H$. It is therefore crucial to choose a small value for $\gamma_H$. The chosen value $\gamma_H = 10^{-3}$ yields an appropriate reduction in spurious velocities with no negative effects on the benchmark values. \\
Figure~\ref{fig: extended velocity} shows the ghost penalty extensions of the velocity and the pressure from the physical domains to the fictitious domains for $t=1.225$. We observe that the ghost penalties provide a rich and stable extension and yield a stable solution independent of the cut location of the fluid interface. The extensions and doubling of cells in the interface layer enable us to capture sharp pressure jumps and velocity kinks at the interface (see Figure~\ref{fig: pressure jump}). 
\paragraph{Test case 2} Figure~\ref{fig: rising bubble results case 2} shows the center of mass, rise velocity and circularity with mesh refinement in comparison to FreeLIFE for mesh size $h=1/160$. We achieve excellent agreement of our results with FreeLIFE until $t \approx 1.5$, which coincides with the onset of the liquid strand formation in the wake of the bubble (see Figure~\ref{fig: rising bubble results case 2 shape}). Similar discrepancies between the results were observed in the original benchmark \cite{Hysing2009} and no agreement for the formation and size of the liquid filaments could be reached in \cite{Hysing2009}. This is due to the fact that the liquid strands breaking or not breaking is highly dependent on the mesh resolution. In our case once thin liquid filaments are contained within one element this part of the strand is no longer detected and will cause the strand to break. The breaking of the liquid strands causes oscillations in the circularity values (see Figure~\ref{fig: rising bubble results case 2} (c)  for $t>2.5$).  For $h=1/40$ and $h=1/80$, the liquid strands break while for $h=1/160$ they do not break (see Figure~\ref{fig: rising bubble results case 2 shape}(a)). Figure~\ref{fig: rising bubble results case 2 shape}(b) shows the influence of the  velocity stabilisation parameters $\gamma_u=\gamma_{div}=\gamma_{gu}=\{0.001,0.01,1\}$ onto the liquid filaments for the final time $t=3.0$ and $h=1/80$. We observe that increasing the velocity stabilisation parameters leads to a thickening of the liquid strands and that no break up occurs for $\gamma_u=1$ for $h=1/80$. This can be explained by interpreting the gradient jump penalties, $j_u, j_{div}$ and $g_u$ as a small amount of additionally viscosity which has this stabilising effect onto the liquid filaments. Note that while an increase in the velocity parameters changes the thickness of the strands, it does not significantly impact the bulk part of the bubble (see Figure~\ref{fig: rising bubble results case 2 shape}(b)).

\begin{table}
\begin{center}
\begin{tabular}{c c c c c c}    \toprule
Case & $\gamma$ & $\rho_1$ & $\rho_2$ & $\eta_1$ & $\eta_2$  \\\midrule
1   & $24.5$  & 100  & 1000  & 1 & 10 \\ 
\rowcolor{black!20} 2 & $1.96$ & 1 & 1000 & 0.1 & 10\\
 \hline
\end{tabular}
\caption{Summary of material parameters for rising bubble benchmark.}
\label{tab: rising bubble cases}
\end{center}
\end{table}

\subsection{Droplet in microfluidic 5:1:5 contraction and expansion}
\begin{figure}
\centering
\includegraphics[width=\textwidth]{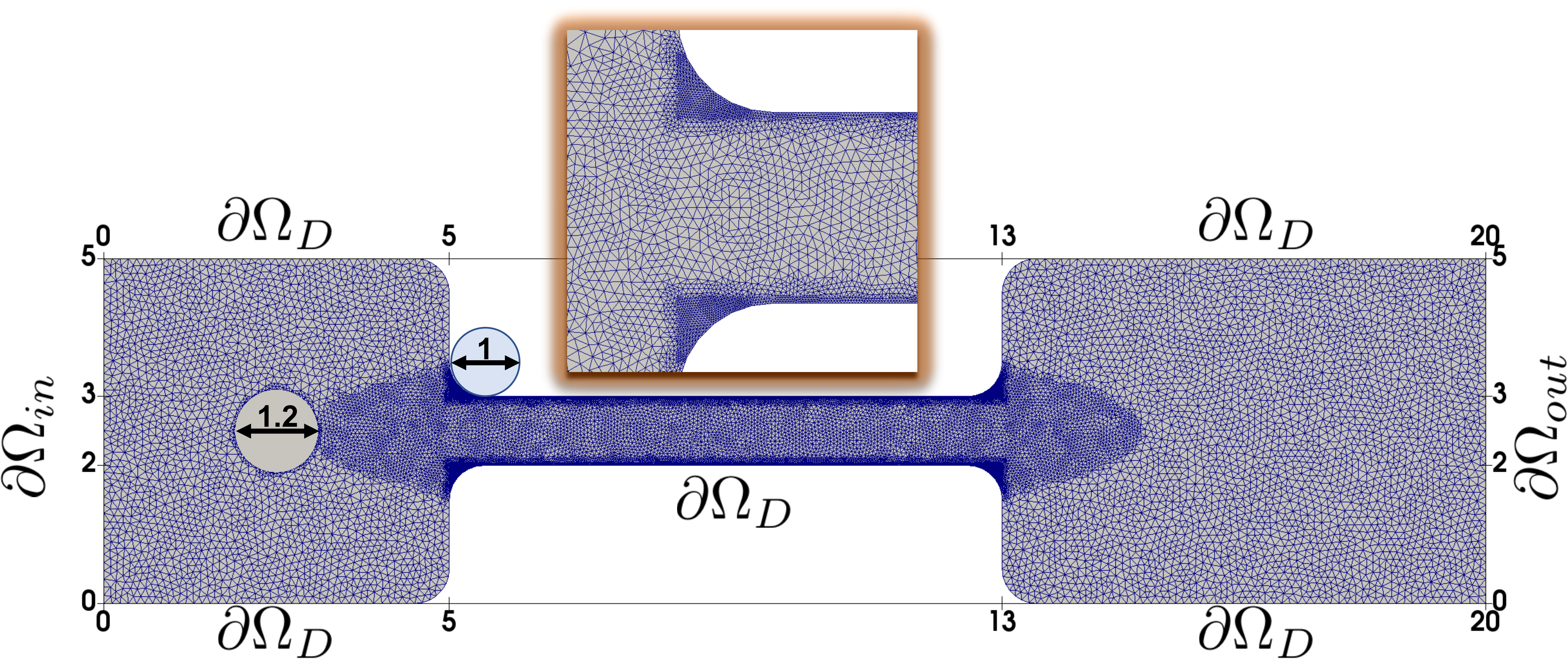}
\caption{Geometry and mesh for the 5:1:5 contraction-expansion problem.}
\label{fig: contraction geometry}
\end{figure}

\begin{figure}
\centering
\subfloat[Case 1, $Ca=0.5$.]{\includegraphics[width=.9\textwidth]{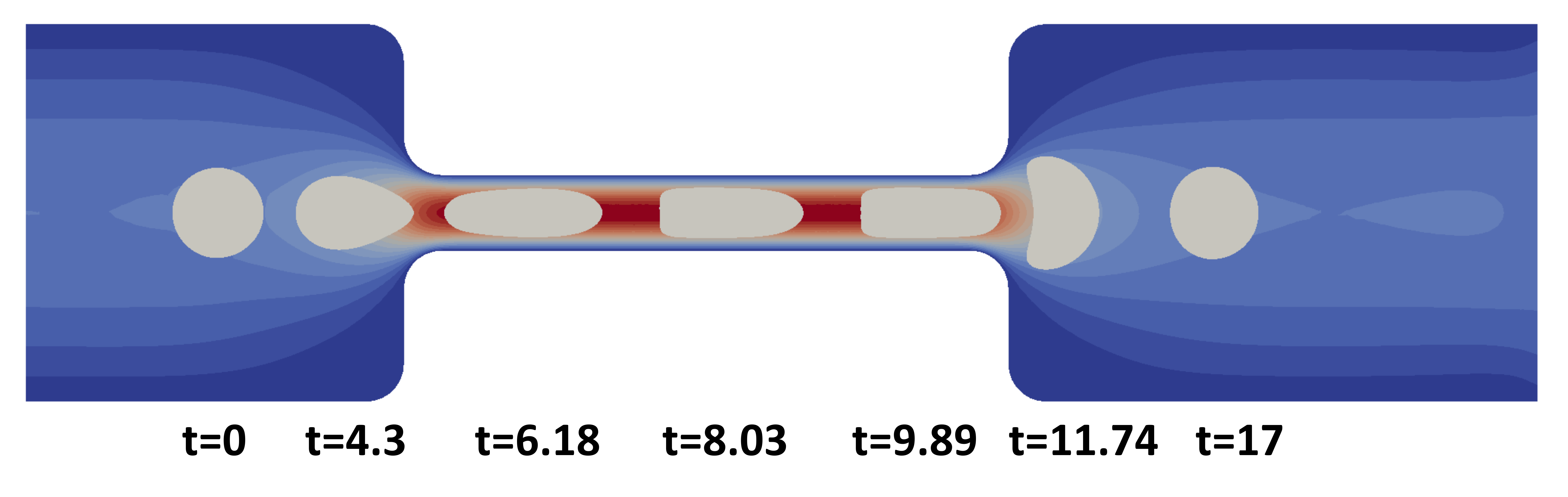}} \\
\subfloat[Case 2, $Ca=0.05$.]{\includegraphics[width=.9\textwidth]{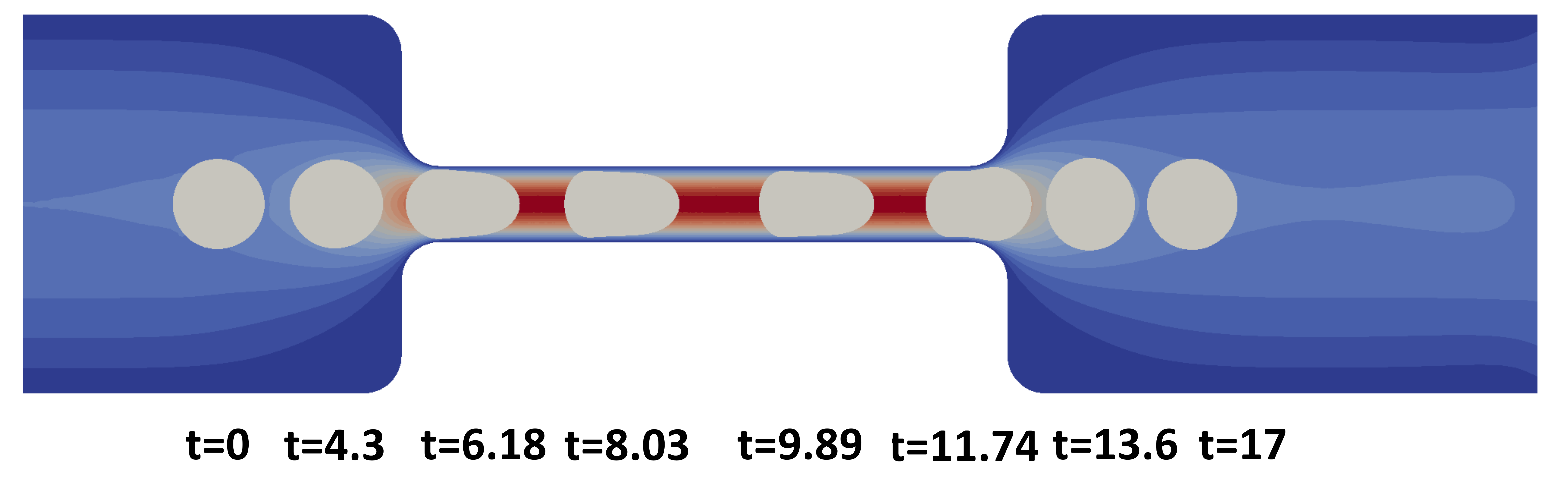}}\\
\subfloat[Case 3, $Ca=0.005$.]{\includegraphics[width=.9\textwidth]{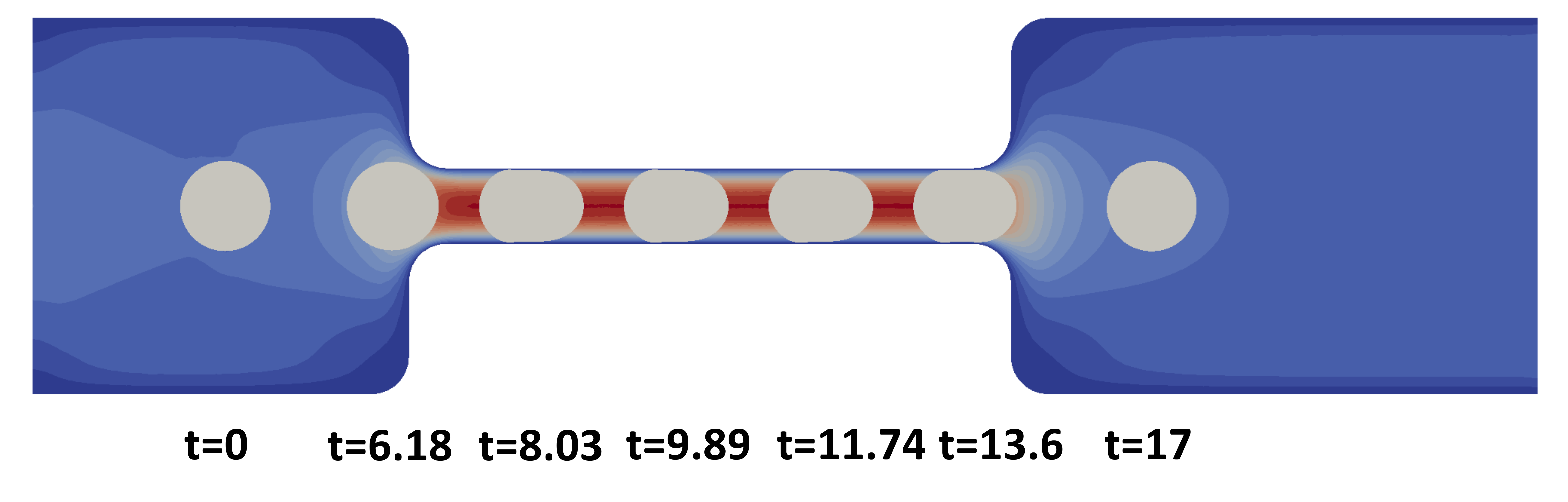}}\\
\subfloat{\includegraphics[width=.4\textwidth]{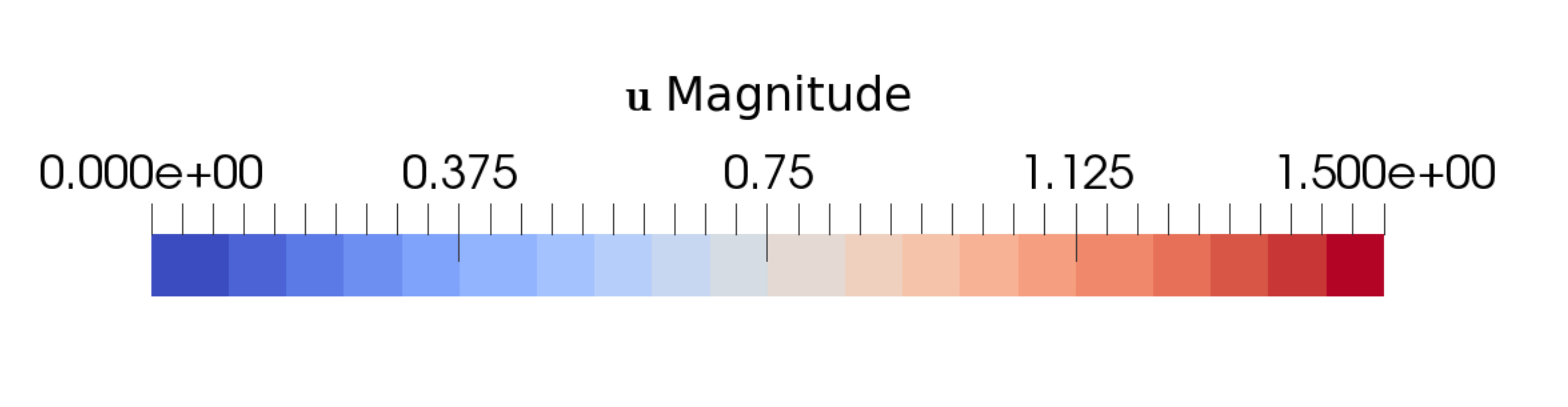}}
\caption{Development of droplet shapes in 5:1:5 contraction-expansion for $Ca=0.5,0.05,0.005$.}
\label{fig: micro-contraction development}
\end{figure}

\begin{figure}
\centering
\subfloat[Case 1, $Ca=0.5$.]{\includegraphics[width=.3\textwidth]{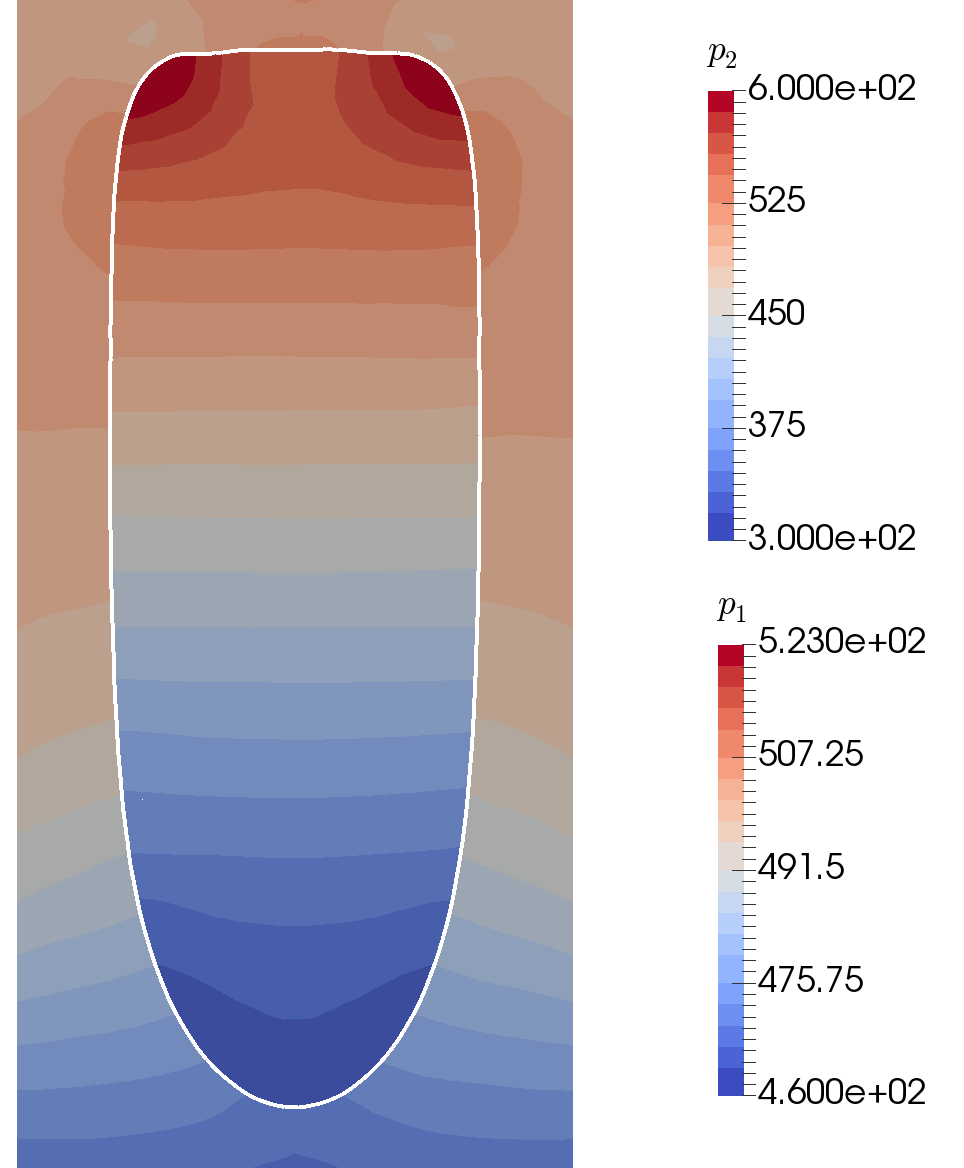}
\includegraphics[width=.3\textwidth]{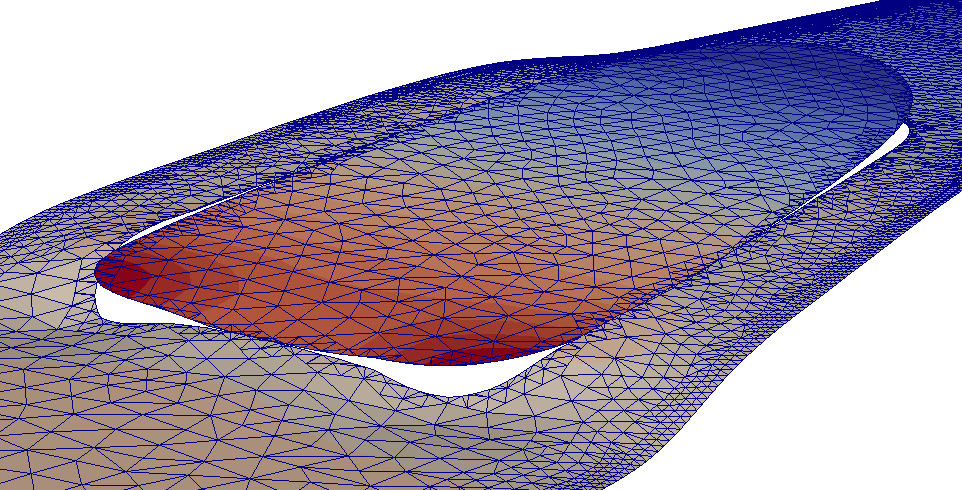}} \\
\subfloat[Case 2, $Ca=0.05$.]{\includegraphics[width=.3\textwidth]{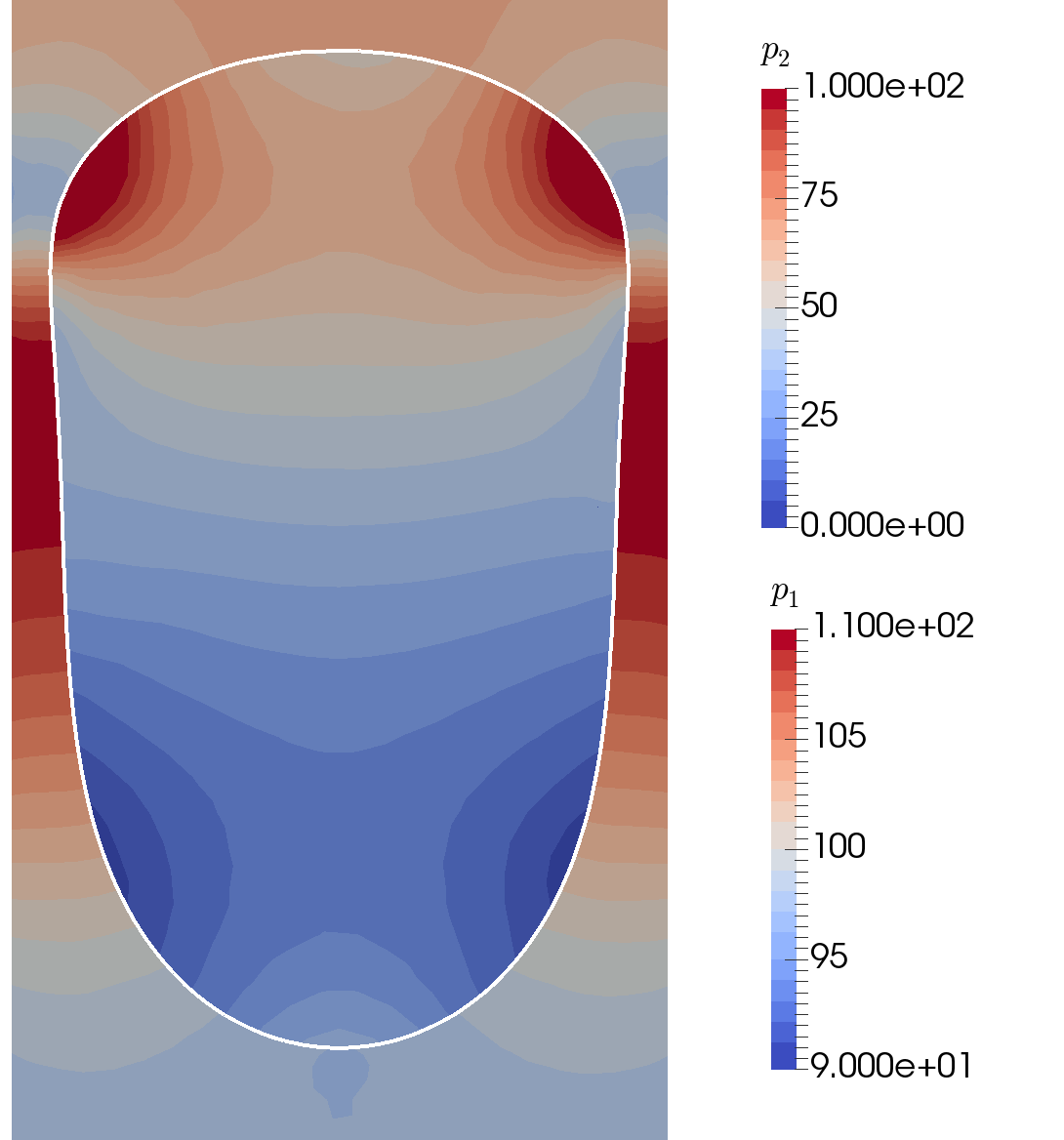}\includegraphics[width=.3\textwidth]{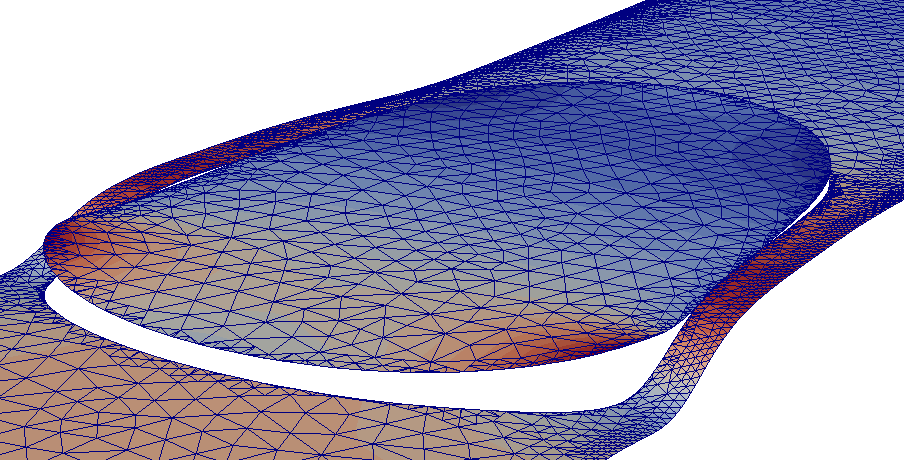}} \\
\subfloat[Case 3, $Ca=0.005$.]{\includegraphics[width=.3\textwidth]{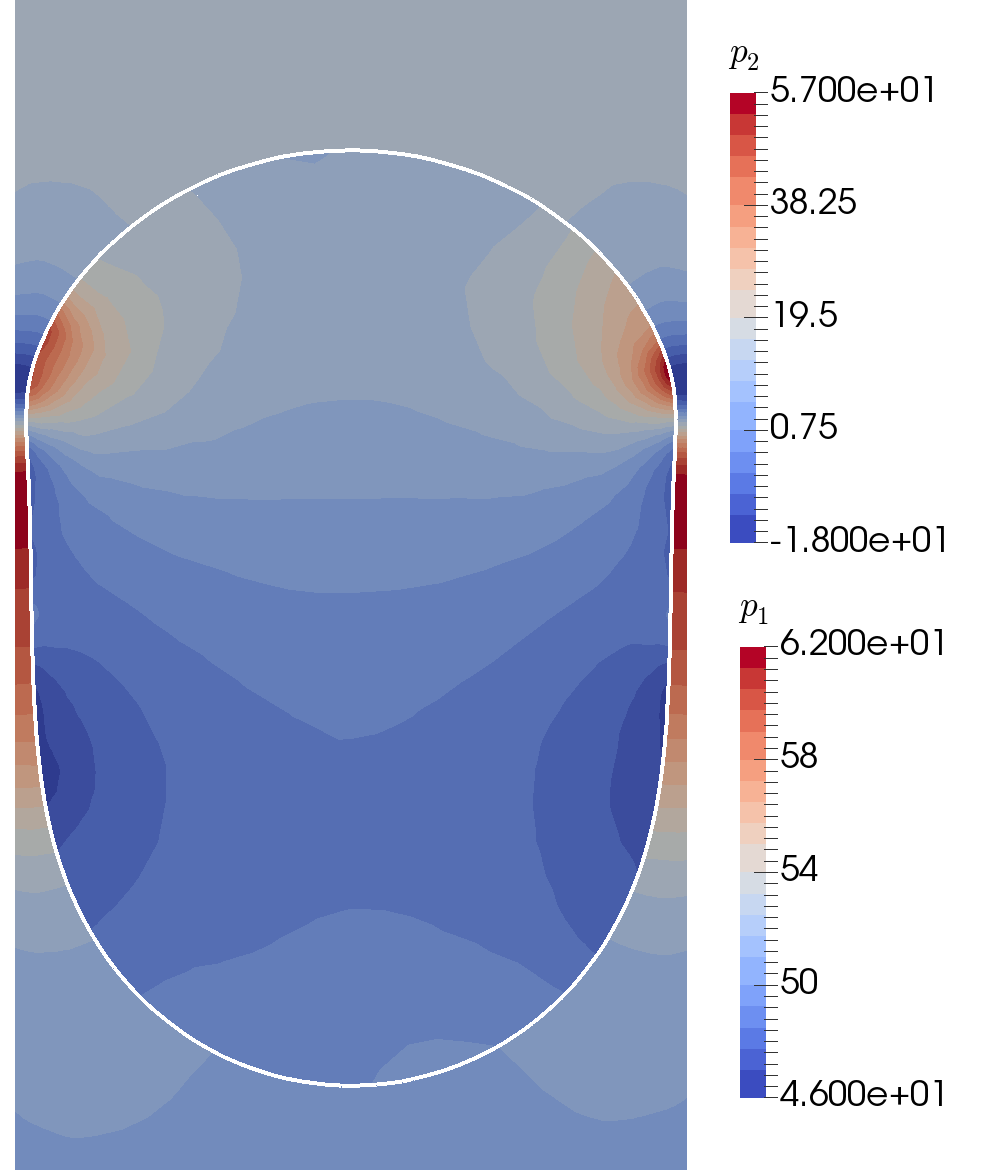}\includegraphics[width=.3\textwidth]{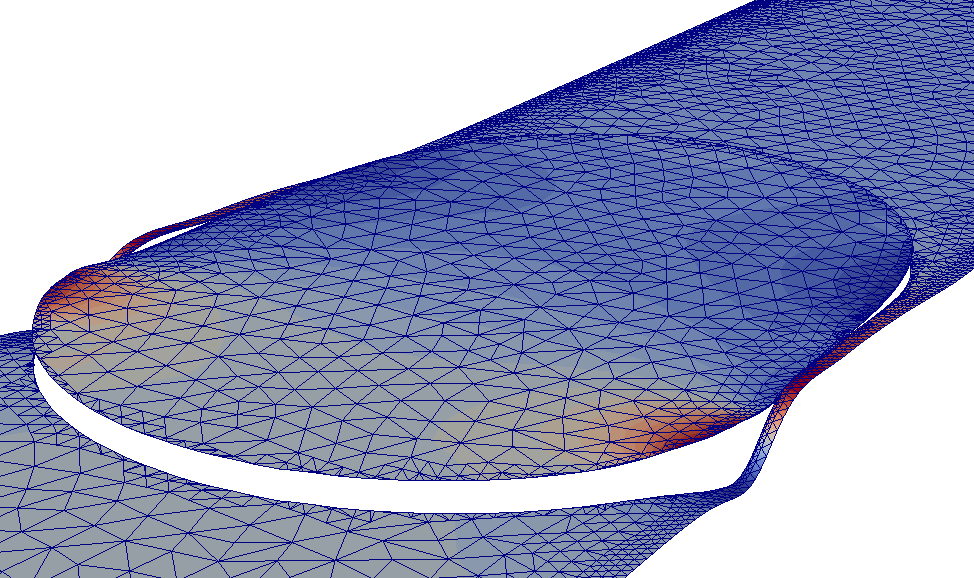}}
\caption{Pressure profiles and droplet shape in middle of narrow channel for  5:1:5 contraction-expansion with decreasing Capillary number.  }
\label{fig: micro-contraction pressure}
\end{figure}

In this section, we consider the deformation of a single droplet in a 5:1:5 contraction and expansion (see e.g. \cite{Chung2009},\cite{Harvie2008},\cite{Hoang2018}) flow. We consider a domain decomposed into an entry domain of height $5$ and length $5$ attached to a narrow channel of height $1$ and length $8$ and an expansion region of height $5$ and length $7$ (see Figure~\ref{fig: contraction geometry}). We consider a droplet placed in the middle of the entry domain $\pmb{x}_m=(2.5,2.5)^T$ of diameter larger than the narrow channel ($d=1.2$). We round the corners of the contraction and expansion with circles of radius $r=0.5$ to prevent stress singularities at these corners. At the top and bottom boundary, we set homogenous Dirichlet boundary condition for the velocity, i.e. $\mathbf{u}=\mathbf{0}$ on $\partial \Omega_D$.  At the left boundary, we set a velocity inflow profile of 
\begin{equation}
\mathbf{u} = (0.048\,(5-y)y, 0)^T \mbox{ on } \partial \Omega_{in}. 
\end{equation}
The scaling factor, $0.048$, is chosen to obtain an average velocity of $\bar{U}=1$. At the right boundary, $\partial \Omega_{out}$, we set homogenous Neumann conditions for the velocity, i.e. $\frac{\partial \mathbf{u}}{\partial \bfn} =0$, and $p=0$. The non-dimensional number, which best characterises this flow configuration is the Capillary number (\cite{Chung2009},\cite{Harvie2008},\cite{Hoang2018}), 
\begin{equation}
Ca = \frac{\eta_2 \bar{U}}{\gamma},
\end{equation}
of the surrounding fluid. The capillary number expresses the ratio of viscous to capillary forces. Here, $L$ is the height of the narrow channel, i.e. $L=1$. As our intend is to simulate conditions in a microfluidic device, where inertial forces are small due to the small dimensions of the geometry, we choose $\rho_1=\rho_2=1.0$. Additionally,  capillary forces tend to dominate viscous forces in microfluidic devices and hence $Ca < 1$. To account for this force balance, we set  the surface tension coefficient to $\gamma=20$ and choose $\eta_2<20$. With the aim of observing different bubble shapes, we fix $\eta_1=1$ and vary $\eta_2$ to create the following three different cases
\begin{eqnarray}
\mbox{Case 1: } \eta_2 =& 10 &\Rightarrow Ca = 0.5, (\eta_1<\eta_2), \\ 
\mbox{Case 2: } \eta_2 =& 1 &\Rightarrow Ca = 0.05,  (\eta_1=\eta_2), \\
\mbox{Case 3: } \eta_2 =& 0.1 &\Rightarrow Ca = 0.005,  (\eta_1>\eta_2).
\end{eqnarray} 
The mesh consisting of a total of 34808 elements is refined once around the narrow channel and refined again towards the channels walls (see Figure~\ref{fig: contraction geometry}). We choose a time step size of $\Delta t = h_{min}/2=0.0062$ and compute the solution in the time interval $t \in[0,17]$. We set the penalty parameters to $\gamma_u=\gamma_{div}=\gamma_{gu}=0.1$, $\gamma_p=0.1$, $\gamma_H = 0.001$, $\hat{\lambda}_{\partial \Omega}=10.0$, $\hat{\lambda}_{\Gamma}=10.0$. \\
Figure~\ref{fig: micro-contraction development} shows the evolution of the droplet through the contraction and expansion regions. The droplet for case 1 undergoes a large deformation and forms a "bullet shape" inside the narrow channel as well as a "half moon shape" in the expansion region. For cases 2 and 3 the droplet takes on a rounder shape inside the channel and relaxes faster to the circular shape in the expansion region. The time interval needed for the droplet to travel through the narrow channel increases from case 1 to case 3, i.e. the mean velocity decreases from case 1 to case 3. Figure~\ref{fig: micro-contraction pressure} shows the pressure profiles and droplet shapes for the droplet in the middle of the narrow channel. We observe a decrease in the fluid layer thickness between the droplet and the wall from case 1 to case 3. The droplet shape becomes increasingly circular from case 1 to case 3. A high pressure forms in the thin liquid layer between the droplet and the channel walls which prevents the droplet from touching the walls. In the rear of the droplet, we have a pressure concentration which sharpens with increasing Capillary number. We have a higher curvature in the front of the droplet than in the rear, which for case 1 yields the bullet shape. This contrast in curvatures decreases with decreasing Capillary number.  Figure~\ref{fig: micro-contraction pressure}  shows that the presented cut finite element method is capable of capturing the jumps in the pressure profiles for low Capillary number without spurious oscillations and is capable of resolving the thin liquid layer between the droplet and the wall. 

\section{Conclusions}
We have presented a cut finite element method for two-phase Navier-Stokes problems with a unified continuous interior penalty stabilisation strategy. Ghost penalties are used to extend the velocity and pressure from physical domains to extended fictitious domains for each phase to guarantee a well defined velocity field of the past time step in the current geometry. In addition, we have employed a smoothened curvature computation to reduce spurious velocity oscillations. \\
We have validated our scheme on two problems: the rising bubble benchmark problem \cite{Hysing2009} and a droplet in a 5:1:5 contraction and expansion microfluidic flow. We have found excellent agreement between our scheme and the rising bubble benchmark. However, we have shown that the curvature stabilisation parameter needs to be chosen with care in order to avoid unphysical shape deformations.  We have demonstrated that the strengths of the presented approach are its capability to represent pressure jumps and velocity kinks, to resolve thin lubrication layers as well as its stability for low Capillary number flows. 

\section*{Acknowledgements}
The authors gratefully acknowledge the financial support provided by the Welsh Government and Higher Education Funding Council for Wales through the S\^{e}r Cymru National Research Network in Advanced Engineering and Materials under grants NRNG06 and NRN102. Susanne Claus would also like to thank Erik Burman and Sarah Zahedi  for numerous inspiring and helpful discussions. 

\bibliographystyle{plainnat}
\bibliography{bibliography}

\end{document}

%% file: notation.tex
\newcommand{\N}{\mathbb{N}}

\newcommand{\mesh}{\mathcal{T}_h}
\newcommand{\avg}[1]{\{ #1 \}}
\newcommand{\mean}[1]{\langle #1 \rangle}

\DefineVerbatimEnvironment{code}{Verbatim}{frame=single,rulecolor=\color{blue}}

\newcommand{\bfg}{\boldsymbol{g}}
\newcommand{\bfv}{\boldsymbol{v}}
\newcommand{\bfV}{\boldsymbol{V}}

\newcommand{\bfn}{\boldsymbol{n}}
\newcommand{\bfx}{\boldsymbol{x}}

\newcommand{\bfbeta}{{\boldsymbol \beta}} 

\newcommand{\mcT}{\mathcal{T}}
\newcommand{\mcV}{\mathcal{V}}

\newcommand{\restr}[2]{ \left. #1 \right|_{#2}}

\newcommand{\jump}[1]{\left\llbracket #1\right\rrbracket}
\newcommand{\norm}[1]{\left\| #1\right\|}

\newcommand{\vel}{\mathbf{u}} 
\newcommand{\stress}{\pmb{\sigma}} 

\newcommand{\velt}{\mathbf{v}}